\documentstyle[12pt,twoside, epsf]{article}

\def\o{\overline}

\def\picill#1by#2(#3){\epsffile{#3}}

\makeatletter
\long\def\UR#1{\leavevmode\setbox\@tempboxa\hbox{#1}\@tempdima\fboxrule
    \advance\@tempdima \fboxsep \advance\@tempdima \dp\@tempboxa
   \hbox{\lower \@tempdima\hbox
  {\vbox{\hrule \@height \fboxrule
          \hbox{  \hskip\fboxsep
          \vbox{\vskip\fboxsep \box\@tempboxa\vskip\fboxsep}\hskip
                 \fboxsep\vrule \@width \fboxrule}%
                  }}}}
\makeatother

\makeatletter
\long\def\LR#1{\leavevmode\setbox\@tempboxa\hbox{#1}\@tempdima\fboxrule
    \advance\@tempdima \fboxsep \advance\@tempdima \dp\@tempboxa
   \hbox{\lower \@tempdima\hbox
  {\vbox{
          \hbox{  \hskip\fboxsep
          \vbox{\vskip\fboxsep \box\@tempboxa\vskip\fboxsep}\hskip
                 \fboxsep\vrule \@width \fboxrule}%
                 \hrule \@height \fboxrule}}}}
\makeatother

\makeatletter
\long\def\UL#1{\leavevmode\setbox\@tempboxa\hbox{#1}\@tempdima\fboxrule
    \advance\@tempdima \fboxsep \advance\@tempdima \dp\@tempboxa
   \hbox{\lower \@tempdima\hbox
  {\vbox{\hrule \@height \fboxrule
          \hbox{\vrule \@width \fboxrule \hskip\fboxsep
          \vbox{\vskip\fboxsep \box\@tempboxa\vskip\fboxsep}\hskip
                 \fboxsep }%
                  }}}}
\makeatother

\makeatletter
\long\def\LL#1{\leavevmode\setbox\@tempboxa\hbox{#1}\@tempdima\fboxrule
    \advance\@tempdima \fboxsep \advance\@tempdima \dp\@tempboxa
   \hbox{\lower \@tempdima\hbox
  {\vbox{
          \hbox{\vrule \@width \fboxrule \hskip\fboxsep
          \vbox{\vskip\fboxsep \box\@tempboxa\vskip\fboxsep}\hskip
                 \fboxsep }%
                 \hrule \@height \fboxrule}}}}
\makeatother

\let \ttorg \tt \def \tt{\ttorg \obeyspaces}

\begin{document}

\date{}

\title{\bf Virtual Knot Theory -- Unsolved Problems}

\author{
Roger Fenn\\
University of Sussex\\
Department of Mathematics\\
Sussex BN1 9RF, England\\
rogerf@sussex.ac.uk \\
and\\
Louis H. Kauffman \\
Department of Mathematics, Statistics and
Computer Science\\
University of Illinois at Chicago\\
851 South Morgan St., Chicago IL 60607-7045, USA\\
kauffman@uic.edu\\
and\\
Vassily O. Manturov \\
Moscow State University\\
Department of Mechanics and Mathematics\\
119992, GSP-2,Leninskie Gory, MSU, Moscow, Russia\\
vassily@manturov.mccme.ru
}

%\author{Roger A. Fenn, Louis H. Kauffman and Vassily O. Manturov}

 \maketitle

 \thispagestyle{empty}

\section{Introduction}
The purpose of this paper is to give an introduction to virtual knot theory and to record a collection of research problems that the authors have
found fascinating. The second section of the paper introduces the theory and discusses problems in that context. The third section is a
list of specific problems.
\bigbreak

We would like to take this opportunity to thank the many people who have, at the time of this writing, worked on the theory of virtual knots and 
links. The ones explicitly mentioned or referenced in this paper are: R. S. Avdeev, V. G. Bardakov, A. Bartholomew, D. Bar-Natan, S.
Budden, S. Carter, H. Dye, R. Fenn, R. Furmaniak, M. Gousssarov, J. Green, D. Hrencecin, D. Jelsovsky, M. Jordan, T. Kadokami, 
N. Kamada, S. Kamada, L. Kauffman, T. Kishino, G. Kuperberg,  S. Lambropoulou, V. O. Manturov, S. Nelson, M. Polyak, D. E. Radford, S.
Satoh, J. Sawollek, M. Saito, W. Schellhorn, D. Silver, V. Turaev, V. V. Vershinin, O. Viro, S. Williams, P. Zinn-Justin and J. B.
Zuber. See
\cite{Avdeev,Bardakov,FB,FBu,CS,CSH,CSW1,Carter,Dror,P,Dye,DKMin,DK,D1,D2,FENN,FJK,FR,FRS1,FRS2,GPV,HR,HRK,KADOKAMI,Kamada,NKamada,
NKamamda1,NSKamada,KANENOBU,KS,KL,KaMa,GEN,VKT,SVKT,DVK,SLK,KD,LF,KM,KiSa,KIS,KUP,MALONG,Ma0,Ma1,MACURVES',MAKH,Ma3,MAPOLY,
MAPOLY',MACURVES,MAVI,Ma7,Ma9,Ma10,Ma11,Ma12,RBook,NELSON,Nelson1,SATOH,SAW,SAW2,Schellhorn,SW,SW1,SW2,SW3,MT,T,TURAEV,Ver,JZ}. We
apologize to anyone who was left out of this list of participant researchers, and we hope that the  problems described herein will
stimulate people on and off this list to enjoy the beauty of virtual knot theory!
\bigbreak

\noindent {\bf Acknowledgement.} 
It gives the second author great pleasure to acknowledge support from NSF Grant DMS-0245588.
\bigbreak

\section{Virtual Knot Theory}

Knot theory
studies the embeddings of curves in three-dimensional space. Virtual knot theory studies the  embeddings of curves in thickened
surfaces of arbitrary genus, up to the addition and removal of empty handles from the surface.  Virtual knots
have a special  diagrammatic theory that makes handling them very similar to the handling of classical knot diagrams. In fact,
this diagrammatic theory simply involves adding a new type of crossing to the knot diagrams, a {\em virtual crossing} that is
neither under nor over. From a combinatorial point of view, the virtual crossings are artifacts of the representation of the
virtual knot or link in the plane. The extension of the Reidemeister moves that takes care of them respects this viewpoint.
A virtual crossing  (See Figure 1) is represented by two crossing arcs with a small circle
placed around the crossing point.
\bigbreak

Moves on virtual diagrams generalize the Reidemeister moves for classical knot and link
diagrams.  See Figure 1.  One can summarize the moves on virtual diagrams by saying that the classical crossings interact with
one another according to the usual Reidemeister moves. One adds the detour moves for consecutive sequences of virtual
crossings and this completes the description of the moves on virtual diagrams. It is a consequence of  moves (B) and (C)
in Figure 1 that an arc going through any consecutive sequence of virtual crossings
can be moved anywhere in the diagram keeping the endpoints fixed;
the places  where the moved arc crosses the diagram become new virtual crossings. This replacement
is the {\em detour move}. See Figure 1.1.
\bigbreak

One can generalize many structures in classical knot theory to the virtual domain, and use the virtual knots to test
the limits of classical problems such as  the question whether the Jones polynomial detects knots and the classical Poincar\'{e}
conjecture. Counterexamples to these conjectures exist in the virtual domain, and it is an open problem whether any of these
counterexamples are equivalent (by addition and subtraction of empty handles) to  classical knots and links. Virtual knot theory
is a significant domain to be investigated for its own sake and for a deeper understanding of classical knot theory.
\bigbreak

Another way to understand the meaning of virtual diagrams is to regard them as representatives for oriented Gauss codes
(Gauss diagrams) \cite{VKT,GPV}. Such codes do not always have planar realizations and an attempt to embed such a code in the plane
leads to the production of the virtual crossings. The detour move makes the particular choice of virtual crossings
irrelevant. Virtual equivalence is the same as the equivalence relation generated on the collection
of oriented Gauss codes modulo an abstract set of Reidemeister moves on the codes.
\bigbreak

One intuition for virtual knot theory is the idea of a particle
moving in three dimensional space in a trajectory that
occasionally disappears, and then reappears elsewhere. By
connecting the disappearance points and the reappearance points
with detour lines in the ambient space we get a picture of the
motion, but the detours, being artificial, must be treated as
subject to replacements. \bigbreak

{\tt    \setlength{\unitlength}{0.92pt}
\begin{picture}(389,343)
\thicklines   \put(77,1){\framebox(252,95){}}
              \put(199,97){\framebox(189,245){}}
              \put(1,97){\framebox(191,240){}}
              \put(197,55){\vector(1,0){19}}
              \put(211,55){\vector(-1,0){17}}
              \put(278,131){\vector(1,0){19}}
              \put(292,131){\vector(-1,0){17}}
              \put(263,237){\vector(1,0){19}}
              \put(277,237){\vector(-1,0){17}}
              \put(265,305){\vector(1,0){19}}
              \put(279,305){\vector(-1,0){17}}
              \put(87,166){\vector(1,0){19}}
              \put(101,166){\vector(-1,0){17}}
              \put(64,240){\vector(1,0){19}}
              \put(78,240){\vector(-1,0){17}}
              \put(66,310){\vector(1,0){19}}
              \put(80,310){\vector(-1,0){17}}
              \put(332,2){\makebox(24,25){$C$}}
              \put(355,71){\makebox(23,26){$B$}}
              \put(31,75){\makebox(20,21){$A$}}
              \put(303,30){\circle{18}}
              \put(262,31){\circle{18}}
              \put(142,69){\circle{18}}
              \put(103,69){\circle{18}}
              \put(282,11){\line(1,1){39}}
              \put(241,52){\line(1,-1){41}}
              \put(121,91){\line(1,-1){40}}
              \put(83,51){\line(1,1){38}}
              \put(286,60){\line(1,1){34}}
              \put(242,12){\line(1,1){35}}
              \put(242,91){\line(1,-1){79}}
              \put(82,91){\line(1,-1){79}}
              \put(82,12){\line(1,1){35}}
              \put(127,57){\line(1,1){34}}
              \put(361,124){\circle{18}}
              \put(322,123){\circle{18}}
              \put(342,143){\circle{18}}
              \put(244,145){\circle{18}}
              \put(263,164){\circle{18}}
              \put(225,164){\circle{18}}
              \put(232,220){\circle{18}}
              \put(232,252){\circle{18}}
              \put(230,308){\circle{18}}
              \put(340,104){\line(1,1){39}}
              \put(299,146){\line(1,-1){41}}
              \put(302,104){\line(1,1){77}}
              \put(204,106){\line(1,1){79}}
              \put(241,183){\line(6,-5){43}}
              \put(205,146){\line(1,1){36}}
              \put(212,236){\line(5,-4){48}}
              \put(261,278){\line(-6,-5){49}}
              \put(251,328){\line(-6,-5){43}}
              \put(209,329){\line(1,-1){40}}
              \put(250,289){\line(0,1){40}}
              \put(296,327){\line(1,0){40}}
              \put(336,326){\line(0,-1){40}}
              \put(336,287){\line(-1,0){40}}
              \put(207,278){\line(1,-1){40}}
              \put(247,236){\line(-1,-1){39}}
              \put(338,236){\line(1,2){20}}
              \put(339,236){\line(1,-2){19}}
              \put(299,276){\line(1,-2){19}}
              \put(318,236){\line(-1,-2){20}}
              \put(204,186){\line(1,-1){79}}
              \put(302,183){\line(1,-1){79}}
              \put(184,147){\line(-1,-1){16}}
              \put(144,107){\line(1,1){16}}
              \put(144,107){\line(-1,1){15}}
              \put(104,147){\line(1,-1){18}}
              \put(106,185){\line(1,-1){79}}
              \put(106,106){\line(1,1){35}}
              \put(150,154){\line(1,1){34}}
              \put(87,149){\line(-1,1){15}}
              \put(46,189){\line(1,-1){17}}
              \put(31,174){\line(1,1){15}}
              \put(7,149){\line(1,1){16}}
              \put(52,155){\line(1,1){34}}
              \put(7,110){\line(1,1){35}}
              \put(7,189){\line(1,-1){79}}
              \put(124,241){\line(-1,-2){20}}
              \put(105,281){\line(1,-2){19}}
              \put(145,241){\line(1,-2){19}}
              \put(144,241){\line(1,2){20}}
              \put(36,220){\line(4,-3){24}}
              \put(11,242){\line(6,-5){16}}
              \put(27,254){\line(-4,-3){16}}
              \put(60,282){\line(-6,-5){24}}
              \put(47,240){\line(-1,-1){39}}
              \put(7,282){\line(1,-1){40}}
              \put(143,292){\line(-1,0){40}}
              \put(143,331){\line(0,-1){40}}
              \put(103,332){\line(1,0){40}}
              \put(20,310){\line(-1,-1){15}}
              \put(46,333){\line(-1,-1){16}}
              \put(46,293){\line(0,1){40}}
              \put(5,333){\line(1,-1){40}}
\end{picture}}

\begin{center}
{ \bf Figure 1 -- Generalized Reidemeister Moves for Virtuals} \end{center}
\vspace{3mm}

\begin{center}
{\tt    \setlength{\unitlength}{0.92pt}
\begin{picture}(482,223)
\thicklines   \put(10,172){\line(1,0){196}}
              \put(64,212){\line(0,-1){81}}
              \put(104,213){\line(0,-1){80}}
              \put(143,213){\line(0,-1){81}}
              \put(85,92){\line(0,-1){79}}
              \put(125,92){\line(0,-1){80}}
              \put(371,90){\line(0,-1){80}}
              \put(334,90){\line(0,-1){79}}
              \put(351,211){\line(0,-1){81}}
              \put(391,210){\line(0,-1){80}}
              \put(311,209){\line(0,-1){81}}
              \put(230,171){\line(1,0){41}}
              \put(271,171){\line(0,-1){121}}
              \put(271,50){\line(1,0){160}}
              \put(431,51){\line(0,1){120}}
              \put(431,171){\line(1,0){41}}
              \put(64,172){\circle{20}}
              \put(105,172){\circle{20}}
              \put(144,172){\circle{20}}
              \put(335,49){\circle{18}}
              \put(371,49){\circle{18}}
              \put(164,109){\vector(1,0){79}}
              \put(243,109){\vector(-1,0){80}}
\end{picture}}

{ \bf Figure 1.1 -- Detour Move}
\end{center}

\subsection{Flat Virtual Knots and Links}
Every classical knot or link diagram can be regarded as a $4$-regular plane graph with extra structure at the
nodes. This extra structure is usually indicated by the over and under crossing conventions that give
instructions for constructing an embedding of the link in three dimensional space from the diagram.  If we take the diagram
without this extra structure, it is the shadow of some link in three dimensional space, but the weaving of that link is not
specified. It is well known that if one is allowed to apply the Reidemeister moves to such a shadow (without regard to the types
of crossing since they are not specified) then the shadow can be reduced to a disjoint union of circles. This reduction is
no longer true for virtual links. More precisely, let a {\em flat virtual diagram} be a diagram with virtual crossings as we have
described them and {\em flat crossings} consisting in undecorated nodes of the $4$-regular plane graph. Virtual crossings
are flat crossings that have been decorated by a small circle. Two flat virtual diagrams are {\em equivalent} if there is a
sequence of generalized flat Reidemeister moves (as illustrated in Figure 1) taking one to the other. A generalized
flat Reidemeister move is any move as shown in Figure 1, but one can ignore the over or under crossing structure.
Note that in studying flat virtuals the rules for changing virtual crossings among themselves and the rules for changing
flat crossings among themselves are identical. However, detour moves as in Figure 1C are available for virtual crossings
with respect to flat crossings and not the other way around.
\bigbreak

We shall say that a virtual diagram {\em overlies} a flat diagram if the virtual diagram is obtained from the flat diagram by
choosing a crossing type for each flat crossing in the virtual diagram. To each virtual diagram $K$ there is an associated
flat diagram $F(K)$ that is obtained by forgetting the extra structure at the classical crossings in $K.$ Note that if $K$
is equivalent to $K'$ as virtual diagrams, then $F(K)$ is equivalent to $F(K')$ as flat virtual diagrams. Thus, if we can
show that $F(K)$ is not reducible to a disjoint union of circles, then it will follow that $K$ is a non-trivial virtual link.
\bigbreak

\begin{center}
{\tt    \setlength{\unitlength}{0.92pt}
\begin{picture}(507,286)
\thicklines   \put(364,1){\makebox(22,23){K}}
              \put(84,1){\makebox(24,23){D}}
              \put(73,242){\makebox(20,20){H}}
              \put(502,160){\line(-1,-1){38}}
              \put(383,40){\line(1,1){70}}
              \put(438,42){\line(-4,5){51}}
              \put(343,160){\line(5,-6){34}}
              \put(300,41){\line(0,1){34}}
              \put(301,180){\line(0,-1){93}}
              \put(262,160){\line(1,0){80}}
              \put(420,180){\line(-1,0){120}}
              \put(299,40){\line(1,0){82}}
              \put(501,200){\line(-1,0){120}}
              \put(381,200){\line(0,-1){120}}
              \put(382,80){\line(-1,0){122}}
              \put(261,82){\line(0,1){78}}
              \put(440,42){\line(1,0){65}}
              \put(421,179){\line(3,-5){83}}
              \put(301,158){\circle{20}}
              \put(382,178){\circle{20}}
              \put(502,200){\line(0,-1){40}}
              \put(416,71){\circle{20}}
              \put(158,73){\circle{20}}
              \put(125,42){\line(1,1){120}}
              \put(244,202){\line(0,-1){40}}
              \put(124,180){\circle{20}}
              \put(43,160){\circle{20}}
              \put(86,162){\line(4,-5){95}}
              \put(163,181){\line(3,-5){83}}
              \put(182,44){\line(1,0){65}}
              \put(3,84){\line(0,1){78}}
              \put(124,82){\line(-1,0){122}}
              \put(123,202){\line(0,-1){120}}
              \put(243,202){\line(-1,0){120}}
              \put(41,42){\line(1,0){82}}
              \put(42,180){\line(0,-1){138}}
              \put(162,182){\line(-1,0){120}}
              \put(4,162){\line(1,0){80}}
              \put(25,241){\oval(16,16)}
              \put(64,264){\line(-1,0){39}}
              \put(64,225){\line(0,1){39}}
              \put(24,225){\line(1,0){40}}
              \put(24,263){\line(0,-1){38}}
              \put(4,283){\line(0,-1){41}}
              \put(43,242){\line(-1,0){39}}
              \put(44,283){\line(0,-1){41}}
              \put(4,283){\line(1,0){40}}
\end{picture}}

{\bf Figure 2 -- Flats $H$ and $D$, and the knot $K.$} \end{center}
\bigbreak

Figure 2 illustrates an example of a flat virtual link $H.$ This link cannot be undone in the flat category because it has an
odd number of virtual crossings between its two components and each generalized Reidemeister move preserves the parity of the
number of virtual crossings between components.  Also illustrated in Figure 2 is a flat diagram $D$ and a virtual knot $K$ that
overlies it. This example is given in \cite{VKT}. The knot shown is undetectable by many invariants (fundamental group, Jones polynomial)
but it is knotted and this can be seen either by using a generalization of the Alexander polynomial that we describe below,
or by showing that the underlying diagram $D$ is a non-trivial flat virtual knot using the filamentation invariant that
is introduced in \cite{HR}. The filamentation invariant is a combinatorial method that is sometimes
successful in indentifying irreducible flat virtuals. At this writing we know very few invariants of flat virtuals.
The flat virtual diagrams present a strong challenge for the construction
of new invariants.  It is important to understand the structure of flat virtual knots and links. This structure
lies at the heart of the comparison of classical and virtual links. We wish to be able to determine when a given virtual
link is equivalent to a classcal link. The reducibility or irreducibility of the underlying flat diagram is the first
obstruction to such an equivalence.
\bigbreak

\subsection{Interpretation of Virtuals as Stable Classes of Links in  Thickened Surfaces}
There is a useful topological interpretation for this virtual
theory in terms of embeddings of links in thickened surfaces. See
\cite{VKT,DVK,KUP}.  Regard each virtual crossing as a shorthand
for a detour of one of the arcs in the crossing through a 1-handle
that has been attached to the 2-sphere of the original diagram. By
interpreting each virtual crossing in this way, we obtain an
embedding of a collection of circles into a thickened surface
$S_{g} \times R$ where $g$ is the number of virtual crossings in
the original diagram $L$, $S_{g}$ is a compact oriented surface of
genus $g$ and $R$ denotes the real line.  We say that two such
surface embeddings are {\em stably equivalent} if one can be
obtained from another by isotopy in the thickened surfaces,
homeomorphisms of the surfaces and the addition or subtraction of
empty handles.  Then we have the \smallbreak \noindent {\bf
Theorem \cite{VKT,DKT,KUP}.} {\em Two virtual link diagrams are
equivalent if and only if their correspondent surface embeddings
are stably equivalent.} \smallbreak \noindent \bigbreak

\noindent{\bf Virtual knots and links give rise to a host of
problems.} As we saw in the previous section, there are
non-trivial virtual knots with unit Jones polynomial.  Moreover,
there are non-trivial virtual knots with integer fundamental group
and trivial Jones polynomial. (Fundamental group is defined
combinatorially by generalizing the Wirtinger presentation.) These
phenomena underline the question of how planarity is involved in
the way the Jones polynomial appears to detect classical knots,
and that the relationship of the fundamental group (and peripheral
system) is a much deeper one than the surface combinatorics for
classical knots. It is possible to take the connected sum of two
trivial virtual diagrams and obtain a non-trivial virtual knot
(the Kishino knot).

Here long knots (or, equivalently $1-1$ tangles) come into play.
Having a knot, we can break it at some point and take its ends to
infinity (say, in a way that they coincide with the horizontal
axis line in the plane). One can study isotopy classes of such
knots. A well-known theorem says that in the classical case, knot
theory coincides with long knot theory. However, this is not the
case for virtual knots. By breaking the same virtual knot at
different points, one can obtain non-isotopic long knots \cite{FJK}.
Furthermore, even if the initial knot is trivial, the resulting
long knot may not be trivial. The ``connected sum'' of two
trivial virtual diagrams may not be trivial in the compact case.
The phenomenon occurs because these two knot diagrams may  be
non-trivial in the long category. It is sometimes more convenient
to consider long virtual knots rather than compact virtual knots,
since connected sum is well-defined for long knots. It is
important to construct long virtual knot invariants to see whether
long knots are trivial and whether they are
classical. One approach is to regard long knots as $1-1$ tangles
and use extensions of standard invariants  (fundamental group,
quandle, biquandle, etc). Another approach is to distinguish two
types of crossings: those having early undercrossing and those
having later undercrossing with respect to the orientation of the
long knot. The latter technique is described in \cite{MALONG}.
\bigbreak

Unlike classical knots, the connected sum of long knots is not commutative \cite{RBook,MaLong}.
Thus, if we show that two long knots $K_{1}$ and $K_{2}$ do not
commute, then we see that they are different and both non-classical.

A typical example of such knots is the two parts of the Kishino
knot, see Figure 2.1.

\begin{center}
{\tt    \setlength{\unitlength}{0.92pt}
\begin{picture}(195,239)
\thicklines   \put(111,11){\line(0,1){47}}
              \put(114,99){\line(0,1){47}}
              \put(95,10){\line(0,1){47}}
              \put(92,99){\line(0,1){47}}
              \put(111,58){\line(1,0){29}}
              \put(114,98){\line(1,0){29}}
              \put(65,58){\line(1,0){29}}
              \put(62,99){\line(1,0){29}}
              \put(24,59){\line(0,-1){12}}
              \put(163,79){\circle{20}}
              \put(45,78){\circle{20}}
              \put(24,60){\line(1,1){39}}
              \put(183,45){\line(0,1){17}}
              \put(133,43){\line(1,0){50}}
              \put(133,51){\line(0,-1){8}}
              \put(133,95){\line(0,-1){34}}
              \put(133,114){\line(0,-1){10}}
              \put(184,114){\line(-1,0){49}}
              \put(183,100){\line(0,1){13}}
              \put(77,46){\line(0,1){7}}
              \put(24,45){\line(1,0){53}}
              \put(77,95){\line(0,-1){33}}
              \put(77,113){\line(0,-1){10}}
              \put(25,113){\line(1,0){52}}
              \put(25,99){\line(0,1){13}}
              \put(141,59){\line(1,1){40}}
              \put(145,98){\line(1,-1){37}}
              \put(24,99){\line(1,-1){40}}
              \put(12,212){\line(1,-1){40}}
              \put(53,213){\line(1,0){41}}
              \put(96,211){\line(1,-1){37}}
              \put(52,171){\line(1,0){40}}
              \put(92,172){\line(1,1){40}}
              \put(13,212){\line(0,1){13}}
              \put(13,226){\line(1,0){52}}
              \put(65,226){\line(0,-1){10}}
              \put(65,208){\line(0,-1){33}}
              \put(12,171){\line(0,-1){12}}
              \put(12,158){\line(1,0){53}}
              \put(65,159){\line(0,1){7}}
              \put(134,213){\line(0,1){13}}
              \put(135,227){\line(-1,0){49}}
              \put(84,227){\line(0,-1){10}}
              \put(84,208){\line(0,-1){34}}
              \put(84,164){\line(0,-1){8}}
              \put(84,156){\line(1,0){50}}
              \put(134,158){\line(0,1){17}}
              \put(12,173){\line(1,1){39}}
              \put(33,191){\circle{20}}
              \put(114,192){\circle{20}}
              \put(139,175){\makebox(32,31){$K$}}
\end{picture}}

{\bf Figure 2.1 -- Kishino and Parts}
\end{center}

We have a natural map

$$\langle\mbox{Long virtual knots }\rangle\to \langle\mbox{Oriented compact virtual
knots}\rangle,$$

\noindent obtained by taking two infinite ends of the long knots together to
make a compact knot. This map is obviously well defined.

It map allows one to construct (weak) long virtual knot invariants
from classical invariants, i.e., just to regard compact knot
invariants as long knot invariants. There is no well-defined
inverse for this map. But, if we were able to construct the map
from compact virtual knots to long virtual knots, we could apply
the long techniques for the compact case. This map does have an
inverse for \underline{classical} knots. Thus, the long techniques
are applicable to classical (long) knots. It would be interesting
to obtain new classical invariants from it. The long category can
also be applied for the case of flat virtuals, where all problems
formulated above occur as well. \bigbreak

There are examples of virtual knots that are very difficult to
prove knotted, and there are infinitely many flat virtual diagrams
that appear to be irreducible, but we have no techniques to prove
it. How can one tell whether a virtual knot is classical? One can ask: Are
there non-trivial virtual knots whose connected sum is trivial?
The latter question cannot be shown by classical techniques, but it can be analyzed by using the surface interpretation
for virtuals. See \cite{MAN}.
\bigbreak

In respect to virtual knots, we are in the same
position as the compilers of the original knot tables. We are, in fact, in developing tables.
At Sussex, tables of virtual knots are being constructed, and tables will appear in a book being written by Kauffman and
Manturov.  See also the website `Knotilus" \cite{knotilus} where there is a tablulation of virtual knots initiated by 
Ralph Furmaniak and Louis Kauffman and the `Knot Atlas" of Dror Bar-Natan, containing a subatlas of virtual knots worked out
in collaboration of Dror and Jeremy Green \cite{Dror}. The theory
of invariants of virtual knots, needs more development. Flat virtuals (whose
study is a generalization of the classification of immersions) are
a nearly unknown territory (but see \cite{HR,TURAEV}). The flat virtuals
provide the deepest challenge since we have very few invariants to
detect them. Curiously, there are many invariants of long flat virtual knots, due to the fact that the virtual (long) knot class of 
the descending virtual diagram associated with a long flat is an invariant of the long flat. (This observation is due to Turaev.)
\bigbreak

\section{Jones Polynomial of Virtual Knots}

We use a generalization of the bracket state summation model for the Jones polynomial to extend it to virtual knots
and links.  We call a diagram in the plane
{\em purely virtual} if the only crossings in the diagram are virtual crossings. Each purely virtual diagram is equivalent by the
virtual moves to a disjoint collection of circles in the plane.
\bigbreak

Given a link diagram $K$, a state $S$ of this diagram is obtained by
choosing a smoothing for each crossing in the diagram and labelling that smoothing with either $A$ or $A^{-1}$
according to the convention that a counterclockwise rotation of the overcrossing line sweeps two
regions labelled $A$, and that a smoothing that connects the $A$ regions is labelled by the letter $A$. Then, given
a state $S$, one has the evaluation $<K|S>$ equal to the product of the labels at the smoothings, and one has the
evaluation $||S||$ equal to the number of loops in the state (the smoothings produce purely virtual diagrams).  One then has
the formula
$$<K> = \Sigma_{S}<K|S>d^{||S||-1}$$
where the summation runs over the states $S$ of the diagram $K$, and $d = -A^{2} - A^{-2}.$
This state summation is invariant under all classical and virtual moves except the first Reidemeister move.
The bracket polynomial is normalized to an
invariant $f_{K}(A)$ of all the moves by the formula  $f_{K}(A) = (-A^{3})^{-w(K)}<K>$ where $w(K)$ is the
writhe of the (now) oriented diagram $K$. The writhe is the sum of the orientation signs ($\pm 1)$ of the
crossings of the diagram. The Jones polynomial, $V_{K}(t)$ is given in terms of this model by the formula
$$V_{K}(t) = f_{K}(t^{-1/4}).$$
\noindent The reader should note that this definition is a direct generalization to the virtual category of the
state sum model for the original Jones polynomial. It is straightforward to verify the invariances stated above.
In this way one has the Jones polynomial for virtual knots and links.
\bigbreak

In terms of the interpretation of virtual knots as stabilized classes of embeddings of circles into thickened surfaces,
our definition coincides with the simplest version of the Jones polynomial for links in thickened surfaces. In that
version one counts all the loops in a state the same way, with no regard for their isotopy class in the surface.
It is this equal treatment that makes the invariance under handle stabilization work. With this generalized version of the
Jones polynomial, one has again the problem of finding a geometric/topological interpretation of this invariant. There is
no fully satisfactory topological interpretation of the original Jones polynomial and the problem is inherited by this
generalization.
\bigbreak

\noindent We have in \cite{DVK} the
\smallbreak
\noindent
{\bf Theorem.} {\em To each non-trivial
classical knot diagram of one component $K$ there is a corresponding  non-trivial virtual knot diagram $Virt(K)$ with unit
Jones polynomial.}
\bigbreak

This Theorem is a key ingredient in the problems involving virtual knots. Here is a sketch of its proof.
The proof uses two invariants of classical knots and links that generalize to arbitrary virtual knots and links.
These invariants are the {\em Jones polynomial} and the {\em involutory quandle} denoted by the notation
$IQ(K)$ for a knot or link $K.$
\bigbreak

Given a
crossing $i$ in a link diagram, we define $s(i)$ to be the result of {\em switching} that crossing so that the undercrossing arc
becomes an overcrossing arc and vice versa. We also define the {\em virtualization}
$v(i)$ of the crossing by the local replacement indicated in Figure 3. In this Figure we illustrate how in the virtualization of
the crossing the  original crossing is replaced by a crossing that is flanked by two virtual crossings.
\bigbreak

Suppose that $K$ is a (virtual or classical) diagram with a classical crossing labeled $i.$  Let $K^{v(i)}$ be the diagram
obtained from $K$ by virtualizing the crossing $i$ while leaving the rest of the diagram just as before. Let $K^{s(i)}$ be
the diagram obtained from $K$ by switching the crossing $i$ while leaving the rest of the diagram just as before. Then it
follows directly from the definition of the Jones polynomial that $$V_{K^{s(i)}}(t) = V_{K^{v(i)}}(t).$$
\noindent As far as the Jones
polynomial is concerned, switching a crossing and virtualizing a crossing look the same.
\bigbreak

The involutory quandle \cite{KNOTS} is an algebraic invariant
equivalent to the fundamental group of the double branched cover of a knot or link in the classical case. In this algebraic
system one associates a generator of the algebra $IQ(K)$ to each arc of the diagram $K$ and there is a relation of the form
$c = ab$ at each crossing, where $ab$ denotes the (non-associative) algebra product of $a$ and $b$ in $IQ(K).$ See Figure 4.
In this Figure we have illustrated through the local relations the fact that  $$IQ(K^{v(i)}) = IQ(K).$$
\noindent As far the involutory quandle is concerned, the original crossing and the virtualized crossing look the same.
\bigbreak

If a classical knot is actually knotted, then its involutory quandle is non-trivial \cite{W}. Hence if we start
with a non-trivial classical knot, we can virtualize any subset of its crossings to obtain a virtual knot that is still
non-trivial. There is a subset $A$ of the crossings of a classical knot $K$ such that the knot $SK$ obtained by
switching these crossings is an unknot.  Let $Virt(K)$ denote the virtual diagram obtained from $A$ by virtualizing
the crossings in the subset $A.$  By the above discussion the Jones polynomial of $Virt(K)$ is the same
as the Jones polynomial of $SK$, and this is $1$ since $SK$ is unknotted. On the other hand, the $IQ$ of $Virt(K)$ is the
same as the $IQ$ of $K$, and hence if $K$ is knotted, then so is $Virt(K).$   We have shown that $Virt(K)$ is a non-trivial
virtual knot with unit Jones polynomial.  This completes the proof of the Theorem.
\bigbreak

{\tt    \setlength{\unitlength}{0.92pt}
\begin{picture}(326,162)
\thinlines    \put(190,141){\line(0,-1){82}}
              \put(242,34){\makebox(41,41){s(i)}}
              \put(170,10){\makebox(41,41){v(i)}}
              \put(1,35){\makebox(41,42){i}}
              \put(292,81){\line(1,0){33}}
              \put(246,81){\line(1,0){35}}
              \put(286,161){\line(0,-1){160}}
              \put(143,81){\circle{16}}
              \put(191,80){\circle{16}}
              \put(190,59){\line(-1,0){25}}
              \put(167,22){\line(0,-1){20}}
              \put(144,22){\line(1,0){22}}
              \put(143,102){\line(0,-1){80}}
              \put(165,103){\line(-1,0){22}}
              \put(165,87){\line(0,1){15}}
              \put(165,59){\line(0,1){17}}
              \put(167,141){\line(1,0){23}}
              \put(166,161){\line(0,-1){20}}
              \put(126,81){\line(1,0){80}}
              \put(46,74){\line(0,-1){73}}
              \put(46,161){\line(0,-1){74}}
              \put(6,81){\line(1,0){80}}
\end{picture}}
\begin{center} {\bf Figure 3 --  Switching and Virtualizing a Crossing} \end{center}
\bigbreak

{\tt    \setlength{\unitlength}{0.92pt}
\begin{picture}(214,164)
\thinlines    \put(197,64){\makebox(16,16){b}}
              \put(136,144){\makebox(22,18){c =}}
              \put(17,143){\makebox(22,18){c =}}
              \put(165,143){\makebox(19,20){ab}}
              \put(113,64){\makebox(16,16){b}}
              \put(67,64){\makebox(16,16){b}}
              \put(43,1){\makebox(17,19){a}}
              \put(47,142){\makebox(19,20){ab}}
              \put(1,63){\makebox(16,16){b}}
              \put(165,1){\makebox(17,19){a}}
              \put(185,142){\line(0,-1){82}}
              \put(138,82){\circle{16}}
              \put(186,81){\circle{16}}
              \put(185,60){\line(-1,0){25}}
              \put(162,23){\line(0,-1){20}}
              \put(139,23){\line(1,0){22}}
              \put(138,103){\line(0,-1){80}}
              \put(160,104){\line(-1,0){22}}
              \put(160,88){\line(0,1){15}}
              \put(160,60){\line(0,1){17}}
              \put(162,142){\line(1,0){23}}
              \put(161,162){\line(0,-1){20}}
              \put(121,82){\line(1,0){80}}
              \put(41,75){\line(0,-1){73}}
              \put(41,162){\line(0,-1){74}}
              \put(1,82){\line(1,0){80}}
\end{picture}}
\begin{center} {\bf Figure 4 --   $IQ(Virt(K) = IQ(K)$} \end{center}
\bigbreak

If there exists a classical knot with unit Jones polynomial, then
one of the knots $Virt(K)$ produced by this theorem may be
equivalent to  a classical knot.  It is an intricate task to
verify that specific examples of $Virt(K)$ are not classical. This
has led to an investigation of new invariants for virtual knots.
In this investigation a number of issues appear. One can examine
the combinatorial generalization of the fundamental group (or
quandle) of the virtual knot and sometimes one can prove by pure
algebra that the resulting group is not classical. This is related
to observations by Silver and Williams \cite{SW}, Manturov
\cite{MAPOLY, MAPOLY'} and by Satoh \cite{SATOH} showing that the
fundamental group of a virtual knot can be interpreted as the
fundamental group of the complement of a torus embedded in four
dimensional Euclidean space. A very fruitful line of new
invariants comes about by examining a generalization of the
fundamental group or quandle that we call the {\em biquandle} of
the virtual knot. The biquandle is discussed in the next Section.
Invariants of flat knots (when one has them) are useful in this
regard. If we can verify that the flat knot $F(Virt(K))$ is
non-trivial, then $Virt(K)$ is non-classical. In this way the
search for classical knots with unit Jones polynomial expands to
the exploration of the structure of the infinite collection of
virtual knots with unit Jones polynomial. \bigbreak

Another way of putting this theorem is as follows: In the arena of
knots in thickened surfaces there are many examples of knots with
unit Jones polynomial. Might one of these be equivalent via handle
stabilization to a classical knot? In \cite{KUP} Kuperberg shows
the uniqueness of the  embedding of minimal genus in the stable
class for a given virtual link. The minimal embedding genus can be
strictly less than the number of virtual crossings in a diagram
for the link.  There are many problems associated with this
phenomenon. \bigbreak

There is a generalization of the Jones polynomial that
involves surface representation of virtual knots. See \cite{Dye,DKMin, MACURVES,MACURVES'}.
These invariants essentially use the fact that the Jones polynomial can be extended to knots in thickened surfaces by keeping
track of the isotopy classes of the loops in the state summation  for this polynomial. In the approach of Dye and Kauffman, one
uses this generalized polynomial directly. In the approach of Manturov, a polynomial invariant is defined using the stabilization
description of the virtual knots.
\bigbreak

\subsection{Atoms}
 An {\em atom} is a pair: $(M^{2},\Gamma)$ where
$M^{2}$ is a closed $2$-manifold and $\Gamma$ is a 4-valent graph
in $M^{2}$ dividing $M^{2}$ into cells such that these cells admit
a checkerboard coloring (the coloring is also fixed). $\Gamma$ is
called the {\em frame} of the atom. See \cite{RBook}, \cite{MAK}.
\smallbreak

Atoms are considered up to natural equivalence, that is, up to 
homeomorphisms of the underlying manifold $M^{2}$ mapping the
frame to the frame and black cells to black cells. 
From this point of view, an atom can be recovered from the frame
together with the following combinatorial structure:

\begin{enumerate}
\item $A$-structure: This indicates which edges for each vertex are
{\em opposite edges}. That is, it indicates the cyclic structure at the vertex.
\item $B$-structure: This indicates pairs of ``black angles''. That is, one
divides the four edges emanating from each vertex into two sets
of adjacent (not opposite) edges such that the black cells are locally
attached along these pairs of adjacent edges.
\end{enumerate}

Given a virtual knot diagram, one can construct the corresponding
atom as follows. Classical crossings correspond to the vertices of
the atom, and generate both the A-structure and the B-structure at these
vertices (the B-structure comes from over/under information). 
Thus, an atom is uniquely determined by a virtual knot diagram.
It is easy to see that the inverse operation is well defined
modulo virtualization. Thus the atom knows everything about
the bracket polynomial (Jones polynomial) of the virtual link.

The crucial notions here are the minimal genus of the atom
and the orientability of the atom.
For instance, for each link diagram with a corresponding
orientable atom (all classical link diagrams are in this class), 
all degrees of the bracket are congruent modulo four
while in the non-orientable case they are congruent only modulo
two.
\smallbreak

The orientability condition is crucial in the construction
of the Khovanov homology theory for virtual links as in \cite{RBook,MAK}.
\bigbreak

\subsection{Biquandles}
In this section we give a sketch of some recent approaches to invariants of virtual knots and links.
\bigbreak

A {\em biquandle} \cite{CS,DVK,FJK,GEN,FB,FBu} is an algebra with 4 binary
operations written $a^b, a_b, a^{\overline{b}}, a_{\overline{b}}$ together
with some relations which we will indicate below. The {\em
fundamental} biquandle is associated with a link diagram and is
invariant under the generalized Reidemeister moves for virtual knots
and links.  The operations in this algebra are motivated by the
formation of labels for the edges of the diagram.  View Figure 5.  In
this Figure we have shown the format for the operations in a
biquandle. The overcrossing arc has two labels, one on each side of
the crossing. There is an algebra element labeling each {\em edge} of
the diagram.  An edge of the diagram corresponds to an edge of the
underlying plane graph of that diagram.  \bigbreak

Let the edges oriented toward a crossing in a diagram be called the {\em input} edges for the crossing, and
the edges oriented away from the crossing be called the {\em output} edges for the crossing.
Let $a$ and $b$ be the input edges for a positive crossing, with $a$ the label of the undercrossing
input and $b$ the label on the overcrossing input. In the biquandle, we label the undercrossing
output by $$c=a^{b},$$ while the overcrossing output is labeled
$$d= b_{a}.$$ \smallbreak

The labelling for the negative crossing is similar using the other two operations.

\noindent To form the fundamental biquandle, $BQ(K)$, we take one generator for each edge of the diagram and two
relations at each crossing (as described above).
\vspace{3mm}

{\tt    \setlength{\unitlength}{0.92pt}
\begin{picture}(280,163)
\thicklines   \put(131,162){\line(0,-1){159}}
              \put(1,1){\framebox(278,161){}}
              \put(231,104){\makebox(21,19){$a^{\o{b}} = a\, \UL{b}$}}
              \put(167,59){\makebox(20,15){$b_{\o{a}} = b\, \LL{a}$}}
              \put(233,62){\makebox(20,16){$b$}}
              \put(212,38){\makebox(18,17){$a$}}
              \put(73,106){\makebox(23,23){$a^{b} = a\, \UR{b}$}}
              \put(77,56){\makebox(32,22){$b_{a} = b\, \LR{a}$}}
              \put(10,86){\makebox(20,16){$b$}}
              \put(52,43){\makebox(18,17){$a$}}
              \put(211,89){\vector(0,1){33}}
              \put(211,43){\vector(0,1){32}}
              \put(250,82){\vector(-1,0){78}}
              \put(51,90){\vector(0,1){33}}
              \put(51,43){\vector(0,1){34}}
              \put(11,83){\vector(1,0){80}}
\end{picture}}

\begin{center}
{\bf Figure 5 -- Biquandle Relations at a Crossing} \end{center}  \vspace{3mm}

\noindent Another way to write this formalism for the biquandle is as follows

$$a^{b} = a\, \UR{b}$$
$$a_{b} = a\, \LR{b}$$
$$a^{\o{b}} = a\, \UL{b}$$
$$a_{\o{b}} = a\, \LL{b}.$$

\noindent We call this the {\em operator formalism} for the biquandle.
\bigbreak

\noindent These considerations lead to the following definition.
\smallbreak

\noindent {\bf Definition.}  A {\em biquandle} $B$ is a set with four binary operations indicated 
above:  $a^{b} \,\mbox{,} \, a^{\o{b}} \, \mbox{,} \,  a_{b} \,\mbox{,} \, a_{\o{b}}.$ We shall refer to the operations with barred
variables as the {\em left} operations and the operations without barred variables as the {\em right} operations. The biquandle is
closed under these operations and the following axioms are satisfied:

\begin{enumerate}
\item Given an element $a$ in $B$, then there exists an $x$ in the biquandle such that $x=a_{x}$ and
$a = x^{a}.$ There also exists a $y$ in the biquandle such that $y=a^{\o{y}}$ and
$a = y_{\o{a}}.$

\item   For any elements $a$ and $b$ in $B$ we have

$$a = a^{b \o{b_{a}}}  \quad \mbox{and} \quad  b= b_{a \o{a^{b}}} \quad \mbox{and}$$

$$a = a^{\o{b}b_{\o{a}}}  \quad \mbox{and} \quad  b= b_{\o{a} a^{\o{b}}}.$$

\item 
\noindent Given elements $a$ and $b$ in $B$ then there exist elements $x, y,
z, t$ such that $x_{b}=a$, $y^{\overline{a}}= b$, $b^x=y$,
$a_{\overline{y}}=x$ and $t^a=b$, $a_t=z$, $z_{\overline{b}}=a$,
$b^{\overline{z}}=t$. The biquandle is called {\em strong} if $x, y,
z, t$ are uniquely defined and we then write $x=a_{b^{-1}},
y=b^{\overline{a}^{-1}}, t=b^{a^{-1}}, z=a_{\overline{b}^{-1}}$,
reflecting the invertive nature of the elements.

\item For any $a$ , $b$ , $c$ in $B$ the following equations hold and the same equations hold when all right operations are
replaced in these equations by left operations.

$$a^{b c} = a^{c_{b} b^{c}} \mbox{,} \quad c_{b a} = c_{a^{b} b_{a}} \mbox{,} \quad (b_{a})^{c_{a^{b}}} = (b^{c})_{a^{c_{b}}}.$$

\end{enumerate}

These axioms are transcriptions of the Reidemeister moves.The first axiom transcribes the first Reidemeister move.
The second axiom transcribes the directly oriented second Reidemeister move.
The third axiom transcribes the reverse oriented Reidemeister move. The fourth axiom transcribes the third Reidemeister move.
Much more work is needed in exploring these  algebras and their applications to knot theory.
\bigbreak

We may simplify the appearance of these conditions by defining
$$S(a,b)=(b_a,a^b),\quad
\overline{S}(a,b)=(b^{\overline{a}},a_{\overline{b}})$$ and in the case
of a strong biquandle,$$S^+_-(a,b)=(b^{a_{b^{-1}}},a_{b^{-1}}),\quad
S^-_+(a,b)=(b^{a^{-1}},a_{b^{a^{-1}}})$$ and 
$${\overline{S}}^{\lower.5ex\hbox{\ $\scriptstyle+$}}_{\ -}(a,b)=
(b_{\ \overline{a^{\overline{b}^{-1}}}}\ , \ a^{\overline{b}^{-1}})=
(b_{\ \overline{a^{b_{a^{-1}}}}} \ , \ a^{b_{a^{-1}}})$$ \ 
and \ $${\overline{S}}^{\lower.5ex\hbox{\ $\scriptstyle-$}}_{\ +}(a,b)=
(b_{\ \overline{a}^{-1}} \ , \ a^{\overline{b_{\overline{a}^{-1}}}})=
(b_{{a}^{b^{-1}}} \ , \ a^{\overline{b_{{a}^{b^{-1}}}}})$$

which we call the {\em sideways} operators. The conditions
then reduce to $$S\overline{S}=\overline{S}S=1,$$ $$ (S\times 1)
(1\times S) (S\times 1) = (1\times S) (S\times 1) (1\times S)$$
$$\overline{S}^-_+S^+_-=S^-_+\overline{S}^+_-=1$$ and finally all
the sideways operators leave the diagonal $$\Delta=\{(a,a)|a\in
X\}$$ invariant.

\subsection{The Alexander Biquandle}

It is not hard to see that
the following equations in a module over $Z[s,s^{-1},t,t^{-1}]$ give a biquandle
structure.

$$a^b=a\,\UR{b} = ta + (1-st)b \, \mbox{,} \quad a_b=a\,\LR{b} = sa$$
$$a^{\overline{b}}=a\,\UL{b} = t^{-1}a + (1-s^{-1}t^{-1})b \, \mbox{,} 
\quad a_{\overline{b}}=a\,\LL{b} = s^{-1}a.$$

\noindent We shall refer to this
structure, with the equations given above, as the {\em Alexander Biquandle}.
\vspace{3mm}

Just as one can define the Alexander Module of a classical knot, we have the Alexander Biquandle of
a virtual knot or link, obtained by taking one generator for each {\em edge} of the projected graph of the  knot diagram and
taking the module relations in the above linear form. Let $ABQ(K)$ denote this module structure for an
oriented link $K$. That is, $ABQ(K)$ is the module generated by the edges of the diagram, factored by the submodule
generated by the relations. This module then has a biquandle structure specified by the operations defined above for an
Alexander Biquandle.
\vspace{3mm}

The determinant of the matrix of relations obtained from the crossings of a diagram
gives a polynomial invariant (up to multiplication by $\pm s^{i}t^{j}$ for integers $i$ and $j$)
of knots and links that we denote by $G_{K}(s,t)$ and call the {\em generalized Alexander polynomial}.
{\bf This polynomial vanishes on classical knots, but is remarkably successful at detecting virtual knots and
links.} In fact $G_{K}(s,t)$ is the same as the polynomial invariant of virtuals
of Sawollek \cite{SAW} and defined by an alternative method by Silver and Williams \cite{SW} and by yet another method by 
Manturov \cite{MAPOLY'}.
It is a reformulation of the invariant for knots in surfaces due to the principal investigator, Jaeger and
Saleur \cite{JKS,KS}.
\vspace{3mm}

We end this discussion of the Alexander Biquandle with two
examples that show clearly its limitations. View Figure 6. In this
Figure we illustrate two diagrams labeled $K$ and $KI.$ It is not
hard to calculate that both $G_{K}(s,t)$ and $G_{KI}(s,t)$ are
equal to zero. However, The Alexander Biquandle of $K$ is
non-trivial -- it is isomorphic to the free module over $Z[s,
s^{-1},t, t^{-1}]$ generated by elements $a$ and $b$ subject to
the relation $(s^{-1} - t -1)(a-b) =0.$ Thus $K$ represents a
non-trivial virtual knot. This shows that it is possible for a
non-trivial virtual diagram to be a connected sum of two trivial
virtual diagrams. However, the diagram $KI$ has a trivial
Alexander Biquandle.  In fact the diagram $KI$, discovered by Kishino \cite{P}, is
now known to be knotted
and its general biquandle is non-trivial. The Kishino diagram has
been shown non-trivial by a calculation of the three-strand Jones
polynomial \cite{KiSa},  by the surface bracket polynomial of Dye and Kauffman \cite{Dye,DKMin}, by the $\Xi$-polynomial (the surface
generalization of the Jones polynomial of Manturov
\cite{MACURVES}, and its biquandle has been shown to be non-trivial by
a quaternionic biquandle representation \cite{FB} of Fenn and Bartholomew which we will now briefly describe. \bigbreak

The quaternionic biquandle is defined by the following operations
where $i^2 = j^2 = k^{2} = ijk = -1, ij = -ji = k,jk = -kj = i, ki
= -ik =j$ in the associative, non-commutative algebra of the
quaternions. The elements $a,b$ are in a module over the
ring of integer quaternions.
$$a \UR{b}=j\cdot a+ (1+i)\cdot b,$$

$$a \LR{b}=-j\cdot a+ (1+i)\cdot b,$$

$$a \UL{b}=j\cdot a+ (1-i) \cdot b,$$

$$a \LL{b}=-j\cdot a+ (1-i)\cdot b.$$
Amazingly, one can verify that these operations satisfy the axioms for the biquandle.
\bigbreak

\noindent Equivalently, referring back to the previous section, define the linear biquandle by
$$S=\pmatrix{
1+i&jt\cr -jt^{-1}&1+i\cr},$$
where $i,j$ have their usual meanings as quaternions and $t$ is a central variable. Let $R$ denote the
ring which they determine. Then as in the Alexander case considered above, for each diagram there is a
square presentation of an $R$-module. We can take the (Study) determinant of the presentation matrix.
In the case of the Kishino knot this is zero. However the greatest common divisor of the codimension 1 
determinants is $2+5t^2+2t^4$ showing that this knot is not classical.

\begin{center} {\tt    \setlength{\unitlength}{0.92pt}
\begin{picture}(264,109)
\thicklines   \put(186,1){\makebox(29,29){$KI$}}
              \put(50,1){\makebox(32,31){$K$}}
              \put(241,70){\circle{20}}
              \put(159,70){\circle{20}}
              \put(105,71){\circle{20}}
              \put(24,70){\circle{20}}
              \put(219,49){\line(-1,0){9}}
              \put(223,89){\line(-1,0){13}}
              \put(180,50){\line(1,0){23}}
              \put(179,90){\line(1,0){24}}
              \put(205,35){\line(1,0){55}}
              \put(206,105){\line(0,-1){70}}
              \put(261,105){\line(-1,0){55}}
              \put(139,90){\line(1,-1){40}}
              \put(223,89){\line(1,-1){37}}
              \put(219,50){\line(1,1){40}}
              \put(140,90){\line(0,1){13}}
              \put(140,104){\line(1,0){52}}
              \put(192,104){\line(0,-1){10}}
              \put(192,86){\line(0,-1){33}}
              \put(139,49){\line(0,-1){12}}
              \put(139,36){\line(1,0){53}}
              \put(192,37){\line(0,1){7}}
              \put(261,91){\line(0,1){13}}
              \put(259,36){\line(0,1){17}}
              \put(139,51){\line(1,1){39}}
              \put(3,52){\line(1,1){39}}
              \put(125,37){\line(0,1){17}}
              \put(75,35){\line(1,0){50}}
              \put(75,43){\line(0,-1){8}}
              \put(75,87){\line(0,-1){34}}
              \put(75,106){\line(0,-1){10}}
              \put(126,106){\line(-1,0){49}}
              \put(125,92){\line(0,1){13}}
              \put(56,38){\line(0,1){7}}
              \put(3,37){\line(1,0){53}}
              \put(3,50){\line(0,-1){12}}
              \put(56,87){\line(0,-1){33}}
              \put(56,105){\line(0,-1){10}}
              \put(4,105){\line(1,0){52}}
              \put(4,91){\line(0,1){13}}
              \put(83,51){\line(1,1){40}}
              \put(43,50){\line(1,0){40}}
              \put(87,90){\line(1,-1){37}}
              \put(44,92){\line(1,0){41}}
              \put(3,91){\line(1,-1){40}}
\end{picture}}
\end{center}

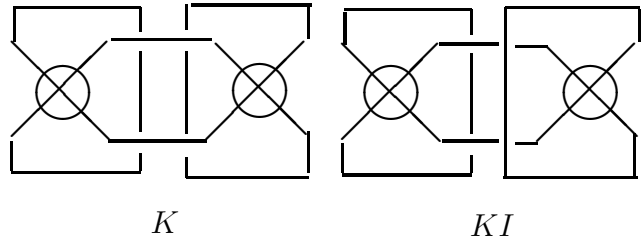
\begin{center} {\bf  Figure 6 -- The Knot $K$ and the Kishino Diagram $KI$ } \end{center}
\vspace{3mm}

\subsubsection{Virtual quandles}

There is another generalization of quandle \cite{Ma1} by means of
which one can obtain the same polynomial as in \cite{SW,SAW} from the
other point of view \cite{MAPOLY,MAPOLY'}. Namely, the formalism is the
same as in the case of quandles at classical crossings but one adds
a special structure at virtual crossings.
The fact that these approaches give the same result in the linear
case was proved recently by Roger Fenn and Andrew Bartholomew.
\bigbreak

Virtual quandles (as well as biquandles) yield generalizations of
the fundamental group and some other invariants. Also, virtual
biquandles admit a generalization for multi-variable polynomials
for the case of multicomponent links, see \cite{MAPOLY}. One can
extend these definitions by bringing together the virtual quandle
(at virtual crossings) and the biquandle (at classical crossings)
to obtain what is called a {\em virtual biquandle}; this work is
now in the process, \cite{KaMa}. \bigbreak

\subsection{A Quantum Model for $G_{K}(s,t),$ Oriented and Bi-oriented Quantum Algebras}

We can understand the structure of the invariant $G_{K}(s,t)$ by rewriting it as a quantum
invariant and then analysing its state summation.  The quantum model for this invariant is obtained in a fashion analogous to the construction of a
quantum model of the Alexander polynomial in \cite{KS,JKS}. The strategy in those papers was to take the basic
two dimensional matrix of the Burau representation, view it as a linear transformation  $T: V \longrightarrow V$ on
a two dimensional module $V$, and them take the induced linear transformation
$\hat{T}: \Lambda^{*}V \longrightarrow  \Lambda^{*}V$ on the exterior algebra of $V$. This gives a transformation on a
four dimensional module that is a solution to the Yang-Baxter equation.  This solution of the Yang-Baxter equation
then becomes the building block for the corresponding quantum invariant.  In the present instance, we have a generalization
of the Burau representation, and this same procedure can be applied to it.
\bigbreak

\noindent The normalized state summation $Z(K)$ obtained by the above process satisfies a skein relation that is just like that of the Conway polynomial:
$Z(K_{+}) - Z(K_{-}) = zZ(K_{0}).$
The basic result behind the correspondence of $G_{K}(s,t)$ and $Z(K)$ is the
 \bigbreak

 \noindent {\bf Theorem \cite{GEN}.} {\em For a (virtual) link $K$, the invariants $Z(K)(\sigma = \sqrt{s}, \tau = 1/\sqrt{t})$ and $G_{K}(s,t)$ are equal
up to
 a multiple of $\pm s^{n}t^{m}$ for integers $n$ and $m$ (this being the well-definedness criterion for $G$).}
 \bigbreak

It is the purpose of this section to place our work with the generalized Alexander polynomial in a context
of bi-oriented quantum algebras and to
introduce the concept of an oriented quantum algebra. In \cite{KRO,KRCAT} Kauffman and
Radford introduce the concept and show that
{\em oriented quantum algebras encapsulate the notion of an oriented quantum link invariant.}
\smallbreak

An {\em oriented quantum algebra} $(A, \rho, D, U)$ is an abstract
model for an oriented quantum invariant of classical links \cite{KRO,KRCAT}. 
For the reader thinking about diagrams, the $\rho$ is an algebraic version of the Yang-Baxter equation and so corresponds to a classical
crossing. The $D$ and $U$ are relatives of cups and caps in the diagrams.
The definition of an oriented quantum algebra is as follows:   We are
given an algebra $A$ over a base ring $k$, an invertible solution
$\rho$ in $A \otimes A$ of the Yang-Baxter equation (in the
algebraic formulation of this equation -- differing from a braiding operator by a transposition), and commuting automorphisms
$U,D:A \longrightarrow A$  of the algebra, such that

$$(U \otimes U)\rho = \rho,$$

$$(D \otimes D)\rho = \rho,$$

$$[(1_{A} \otimes U)\rho)][(D \otimes 1_{A^{op}})\rho^{-1}]
= 1_{A \otimes A^{op}},$$

\noindent and

$$[(D \otimes 1_{A^{op}})\rho^{-1}][(1_{A} \otimes U)\rho)]
= 1_{A \otimes A^{op}}.$$

\noindent The last two equations say that $[(1_{A} \otimes U)\rho)]$ and $[(D \otimes 1_{A^{op}})\rho^{-1}]$
are inverses in the algebra $A \otimes A^{op}$ where $A^{op}$ denotes the opposite algebra.
\vspace{3mm}

When $U=D=T$, then $A$ is said to be {\em balanced}.
In the case where $D$ is the identity mapping, we call the
oriented quantum algebra {\em standard}.  In  \cite{KRCAT} we show that
{\em the invariants defined by Reshetikhin and Turaev (associated with a
quasi-triangular Hopf algebra) arise from standard oriented quantum algebras.}
It is an interesting structural feature of algebras that we have elsewhere
\cite{GAUSS} called {\em quantum algebras} (generalizations of quasi-triangular Hopf algebras)
that they give rise to standard oriented quantum algebras.   It would be of interest to see invariants such as the
Links-Gould invariant \cite{KDL} in this light. 
\bigbreak

\noindent We now extend the concept of oriented quantum algebra  by adding a second solution to
the Yang-Baxter equation $\gamma$ that will take the role of the virtual crossing.
\smallbreak

\noindent {\bf Definition.} A {\em bi-oriented quantum algebra} is a quintuple  $(A, \rho, \gamma, D , U)$ such that
$(A, \rho, D, U)$ and $(A,\gamma, D, U)$ are oriented quantum algebras and $\gamma$ satisfies the following properties:

\begin{enumerate}
\item $\gamma_{12}\gamma_{21} = 1_{A \otimes A}.$ (This is the equivalent to the statement that the
braiding operator corresponding to $\gamma$ is its own inverse.)

\item Mixed identities involving $\rho$ and $\gamma$ are satisfied. These correspond to the
braiding versions of the virtual detour move of type three that involves two virtual crossings and one
real crossing. See \cite{GEN} for the details.

\end{enumerate}
\bigbreak

By extending the methods of \cite{KRCAT}, it is not hard to see that {\em a bi-oriented quantum algebra will always give rise to invariants
of virtual links up to the type one moves (framing and virtual framing).}
\bigbreak

\noindent
In the case of the generalized Alexander polynomial, the state model $Z(K)$ translates directly into a specific example of a bi-oriented balanced quantum
algebra $(A, \rho, \gamma, T).$  The main point about this bi-oriented quantum algebra is that the operator $\gamma$ for the virtual
crossing is {\em not} the identity operator; this non-triviality is crucial to the structure of the invariant.
We will investigate bi-oriented quantum  algebras and other examples of virtual invariants derived from them.
\bigbreak

We have taken a path to explain not only the evolution of a theory of invariants
of virtual knots and links, but also (in this subsection) a description of our oriented quantum algebra formulation of
the whole class of quantum link invariants. Returning to the case of the original Jones polynomial, we want to understand
its capabilities in terms of the oriented quantum algebra that generates the invariant.
\bigbreak

\section{Invariants of Three-Manifolds}
As is well-known, invariants of three manifolds can be formulated in terms of Hopf algebras and quantum algebras and spin
recoupling networks. In formulating such invariants it is useful to represent the three-manifold via surgery on a framed link.
Two framed links that are equivalent in the Kirby calculus of links represent the same three manifold and conversely.
To obtain invariants of three manifolds one constructs invariants of framed links that are also invariant under the Kirby moves
(handle sliding, blowing up and blowing down).

A classical three-manifold is mathematically
the same as a Kirby equivalence class of a framed link. The fundamental group of the three-manifold associated with a link is
equal to the fundamental group of the complement of the link modulo the subgroup generated by the framing longtudes for the link.
We refer to the fundmental group of the three manifold as the {\em three manifold group}.
If there is a counterexample to the
classical Poincar\'{e} conjecture, then the counterexample would be represented by surgery on some link $L$ whose three manifold
group is trivial, but
$L$ is not trivial in Kirby calculus (i.e. it cannot be reduced to nothing).
\bigbreak

Kirby calculus can be generalized to the class of virtual knots and links.
We define a {\em virtual three manifold} to be a Kirby equivalence class of framed virtual links.
The three manifold group generalizes via the
combinatorial fundamental group associated to the virtual link (the framing longitudes still exist for virtual links).
The {\em Virtual Poincar\'{e} Conjecture} to virtuals would say that a virtual three-manifold with trivial
fundamental group is trivial in Kirby calculus. However, {\bf The virtual Poincar\'{e} conjecture is false} \cite{DK}.  There exist
virtual links whose three manifold group is trivial that are nevertheless not Kirby equivalent to nothing. The simplest example
is the virtual knot in Figure 7 .  We detect the non-triviality of the Kirby class of this knot by computing that it has an
$SU(2)$ Witten invariant that is different from the standard three-sphere.
\bigbreak

\begin{center}
{\tt    \setlength{\unitlength}{0.92pt}
\begin{picture}(97,78)
\thicklines   \put(65,18){\line(-1,0){15}}
              \put(65,73){\line(0,-1){55}}
              \put(94,75){\line(-1,0){29}}
              \put(94,55){\line(0,1){20}}
              \put(74,55){\line(1,0){20}}
              \put(22,19){\circle{18}}
              \put(42,43){\line(-1,0){15}}
              \put(4,42){\line(1,0){14}}
              \put(3,20){\line(1,0){36}}
              \put(3,21){\line(0,1){21}}
              \put(21,56){\line(1,0){37}}
              \put(21,3){\line(0,1){53}}
              \put(42,3){\line(-1,0){21}}
              \put(43,43){\line(0,-1){40}}
\end{picture}}
\end{center}

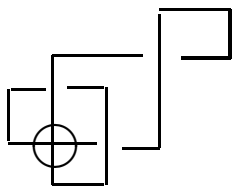
\begin{center} {\bf Figure 7 -- A counterexample to the Poincar\'{e} Conjecture for Virtual Three-Manifolds} \end{center}
\bigbreak

This counterexample to the Poincar\'{e} conjecture in the virtual domain shows how a classical counterexample might behave in the
context of Kirby calculus. Virtual knot theory can be used to search for a counterexample to the classical Poincar\'{e} conjecture by
searching for virtual counterexamples that are equivalent in Kirby calculus to classical knots and links. This is a new and
exciting approach to the dark side of the classical Poincar\'{e} conjecture.
\bigbreak

\section{ Gauss Diagrams and Vassiliev Invariants}
The reader should recall the notion of a {\em Gauss diagram} for a knot. If $K$ is a knot
diagram, then $G(K),$ its Gauss diagram, is a circle comprising the Gauss code of the knot by arranging the traverse of the diagram from
crossing to crossing along the circle and putting an arrow (in the form of a chord of the circle) between the two appearances of the
crossing. The arrow points from the overcrossing segment to the undercrossing segment in the order of the traverse of the diagram. (Note: Turaev uses another convention, \cite{TURAEV}.)
Each chord is endowed with a sign that is equal to the sign of the corresponding crossing in the knot diagram. At the level of the Gauss diagrams,
a virtual crossing is simply the absence of a chord. That is, if we wish to transcribe a virtual knot diagram to a Gauss diagram, we ignore the
virtual crossings. Reidemeister moves on Gauss diagrams are defined by translation from the corresponding diagrams from planar representation.
Virtual knot theory is precisely the theory of {\em arbitrary Gauss diagrams}, up to the Gauss diagram Reidemeister moves. Note that an arbitrary
Gauss diagram is any pattern of directed, signed chords on an oriented circle. When transcribed back into a planar knot diagram, such a Gauss
diagram may require virtual crossings for its depiction.
\bigbreak

In \cite{GPV} Goussarov, Polyak and Viro initiate a very important program for producing Vassiliev invariants of finite type of virtual and classical
knots. The gist of their program is as follows.They define the notion of a semi-virtual crossing, conceived as a dotted, oriented, signed
chord in a Gauss diagram for a knot. An {\em arrow diagram} is a Gauss diagram all of whose chords are dotted. Let $\cal A$ denote
the collection of all linear combinations of arrow diagrams with integer coefficients. Let $\cal G$ denote the collection of all arbitrary Gauss
diagrams (hence all representatives of virtual knots). Define a mapping $$i:\cal G \longrightarrow \cal A$$ by expanding at each chord of
a Gauss diagram $G$ into the sum of replacing the chord by a dotted chord and the removal of that chord. Thus
$$i(G) = \Sigma_{r \in R(G)} G^{r}$$ where $R(G)$ denotes all ways of replacing each chord in $G$ either by a dotted chord, or by nothing; and
$G^{r}$ denotes that particular replacement applied to $G.$
\bigbreak

Now let $\cal P$ denote the quotient of $\cal A$ by the subalgebra generated by the relations in $\cal A$ corresponding to the Reidemeister moves.
Each Reidemeister move is of the form $X = Y$ for certain diagrams, and this translates to the relation $i(X) - i(Y) = 0$ in $\cal P,$ where
$i(X)$ and $i(Y)$ are individually certain linear combinations in $\cal P.$ Let $$I:\cal G \longrightarrow \cal P$$ be the map induced by
$i.$ Then it is a formal fact that $I(G)$ is invariant under each of the Reidemeister moves, and hence that $I(G)$ is an invariant of the
corresponding Gauss diagram (virtual knot) $G.$ The algebra of relations that generate the image of the Reidemeister moves in $\cal P$ is
called the {\em Polyak algebra.}
\bigbreak

So far, we have only desribed a tautological and not a computable invariant. The key to obtaining computable invariants is to truncate.
Let ${\cal P}_{n}$ denote $\cal P$ modulo all arrow diagrams with more than $n$ dotted arrows. Now ${\cal P}_{n}$ is a finitely generated
module over the integers, and the composed map $$I_{n}:{\cal G} \longrightarrow {\cal P}_{n}$$ is also an invariant of virtual knots.
Since we can choose a specific basis for ${\cal P}_{n},$ the invariant $I_{n}$ is in principle computable, and it yields a large collection of
Vassiliev invariants of virtual knots that are of finite type. The paper by Goussarov, Polyak and Viro investigates specific methods for
finding and representing these invariants.  They show that every Vassiliev invariant of finite type for classical knots can be written as
a combinatorial state sum for long knots. They use the virtual knots as an intermediate in the construction.
\bigbreak

By directly constructing Vassiliev invariants of virtual knots from known invariants of virtuals, 
we can construct invariants that are not of
finite type in the above sense (See \cite{VKT}.) These invariants also deserve further investigation.
\bigbreak

\section{Khovanov and Other Invariants}
The {\em Khovanov Categorification of the Jones polynomial} \cite{DB}
is important. This invariant is constructed by
promoting the states in the bracket summation to tensor powers of
a vector space $V$, where a single power of $V$ corresponds to a
single loop in the state. In this way a graded complex is
constructed, whose graded Euler characteristic is equal to the
original Jones polynomial, and the ranks of whose graded homology
groups are themselves invariants of knots. It is now known that
the information in the Khovanov construction exceeds that in the
original Jones polynomial. It is an open problem whether a
Khovanov type construction can generalize to virtual knots in
the general case. The construction for the Khovanov polynomial for virtuals
over ${\bf Z}_{2}$ was proposed in \cite{MAKH}. Khovanov homology (with arbitrary coefficient ring)
can be defined for any link diagram for which the corresponding
atom is orientable. This leads to two explicit geometric constructions
of the virtual link Khovanov homology as in \cite{MAK}. The main point here is that there is a well-formulated theory 
of Khovanov invariants for virtual knots, and it needs more development.
\bigbreak

Recent work by other authors related to knots in thickened surfaces \cite{P} promises to shed light on this
issue. More generally, the subject of upleveling known polynomial--type invariants of knots and links to homology
theories appears to be very fruitful, and new ways to accomplish this in the virtual category should shed light on the nature
of these invariants.
\bigbreak

One of the more promising directions for
relating Vassiliev invariants to our present concerns is the theory of gropes \cite{CT}, where one considers surfaces spanning
a given knot, and then recursively the surfaces spanning curves embedded in the given surface. This hierarchical structure
of curves and surfaces is likely to be a key to understanding the geometric underpinning of the original Jones polynomial.
The same techniques in a new guise could elucidate invariants of virtual knots and links.
\bigbreak

\section{A List of Problems}

Below, we present a list of actual problems
closely connected with virtual knot theory.

\begin{enumerate}

\item {\bf Recognising the Kishino Knot}
There have been invented at least 5 ways to recognize the Kishino
virtual knot (from the unknot): The $3$--strand Jones polynomial
i.e. the Jones polynomial of the $3$--strand cabling of the
knot) \cite{KiSa}, the $\Xi$--polynomial, see \cite{MACURVES}, the quaternionic
biquandle \cite{FB}, and the surface bracket polynomial (Dye and
Kauffman \cite{DKMin}). In \cite{KADOKAMI} Kadokami proves the knot is non-trivial by examining the immersion class of a shadow
curve in genus two.
\bigbreak

Are we done with this knot? Perhaps not. Other proofs of its non-triviality may be illuminating.
The fact that the Kishino diagram is non-trivial and yet a connected sum of trivial virtual knots suggests the question:
{\em Classify when a non-trivial virtual knot can be the connected sum of two trivial virtual knots.}
A key point here is that the connected sum of closed virtual knots is not well defined and hence the different choices give
some interesting effects. With long virtual knots, the connected sum is ordered but well-defined, and the last question is closely
related to the question of classifying the different long virtual knots
whose closures are equivalent to the unknot.
\bigbreak

\item {\bf Flat Virtuals}
See Section $2.1.$ Flat virtual knots, also known as {\em virtual strings} \cite{HR,TURAEV}, are difficult to classify.
Find new combinatorial invariants of flat virtuals.  \smallbreak

We would like to know more about the flat biquandle algebra. This algebra is isomorphic to the Weyl algebra
\cite{FBu} and has no (non-trivial) finite dimensional representations. One can make small examples of the flat biquandle algebra,
that detect some flat linking beyond mod $2$ linking numbers, but the absence of other finite dimensional representations presents
a problem.
\bigbreak

\item {\bf The Flat Hierarchy.} The flat hierarchy is constructed
for any ordinal $\alpha$. We label flat crossings with members of
this ordinal. In a flat third Reidemeister move, a line with two
$a$ labels can slide across a crossing labeled $b$ only if $a$ is
greater than $b.$ This generalizes the usual to theory of flat
virtual diagrams to a system with arbitrarily many different types
of flat crossings. Classify the diagrams in this hierarchy. This concept is due to
L. H. Kauffman (unpublished). A first step in working with the flat hierarchy can be found in 
\cite{Ma11}.

\item {\bf Virtuals and the Theory of Doodles.}
Compare flat theories of virtuals with theories of doodles.
A {\bf doodle link} has only one type of crossing, and that cannot slide over itself.
See \cite{FT,F,KHO}.

\item {\bf Virtual Three Manifolds.} There is a theory of virtual $3$--manifolds constructed as
formal equivalence classes of virtual diagrams modulo generalized
Kirby moves, see \cite{DK}. From this point of view, there are two
equivalences for ordinary $3$--manifolds: homeomorphisms and
virtual equivalence. Do these equivalences coincide? That is,
given two ordinary three manifolds, presented by surgery on framed
links $K$ and $L$, suppose that $K$ and $L$ are equivalent through
the virtual Kirby calculus. Does this imply that they are
equivalent through the classical Kirby calculus?
\bigbreak

\subitem {\em What is a virtual $3$-manifold?} That is, give an
interpretation of these equivalence classes in the domain of
geometric topology. \bigbreak

Construct {\em another} theory of virtual
$3$--manifolds by performing surgery on  links in
thickened surfaces $S_{g}\times R$ considered up
to stabilization. Will this theory coincide to that proposed by Kauffman and Dye
\cite{DK}?

\item {\bf Welded Knots.} We would like to understand welded knots \cite{SATOH}. It is
well known \cite{NELSON,KANENOBU} that if we admit forbidden moves
to the virtual link diagrams, each virtual knot can be transformed
to the unknot. If we allow only one forbidden move (e.g. the upper
one), then there are lots of different equivalence classes of
knots. In fact the fundamental group and the quandle of the
virtual diagram are invariant under the upper forbidden move. The
resulting equivalence classes are called {\em welded knots}. Similarly, {\em  welded braids} were studied
in \cite{FRR}, and every welded knot is the closure of a welded braid. The
question is to construct good invariants of welded knots and, if
possible, to classify them. In \cite{SATOH} a mapping is constructed from welded knots to ambient isotopy classes of
embeddings of tori (ribbon tori to be exact) in four dimensional space, and it is proved that this mapping is an isomorphism from the combinatorial 
fundamental group (in fact the quandle) of the welded knot to the fundamental group of the complement of the corresponding torus embedding
in four-space. Is this Satoh mapping faithful from equivalence classes of welded knots (links) to ambient isotopy classes of ribbon torus embeddings
in four-space?
\bigbreak 

\item {\bf Long Knots and Long Flat Knots.} Enlarge the long knot invariant structure proposed in
\cite{MALONG}. Can one get new classical knot invariants from the approach in this paper?
Bring together the ideas from \cite{MALONG} with the biquandle
construction from \cite{GEN} to obtain more powerful invariants of long
knots. Long flat virtuals can be studied via a powerful remark due to V. Turaev (in conversation) to the effect that 
one can associate to a given long flat virtual knot diagram $F$ a {\em descending} diagram $D(F)$ (by always going over before going
under in resolving the flat (non-virtual) crossings in  the diagram). The long virtual knot type of $D(F)$ is an invariant of the long flat
knot $F.$ This means that one can apply any other invarinant $I$ of virtual knots that one likes to $D(F)$ and $I[D(F)]$ will be an invariant
of the long flat $F.$ It is quite interesting to do sample calculations of such invariants \cite{LF} and this situation underlines the deeper 
problem of finding a full classifiction of long flat knots.
\bigbreak

\item {\bf Virtual Biquandle.} Construct presentations of the virtual biquandle with the a
linear (non-commutative) representation at classical crossings and
some interesting structure at virtual crossings.

\item {\bf Virtual braids.}

%\subitem 
Is there a birack such that its action on virtual
braids is faithful?

%\subitem
Is the invariant of virtual braids in \cite{Ma9}, see also
\cite{Ma0, Bardakov, KL} faithful?

%\subitem
The action defined by {\it linear} biquandles is not faithful. This
almost certainly means that the corresponding linear invariants of
virtual knots and links are not faithful \cite{GQ}.

\item {\bf The Fundamental Biquandle.} Does the fundamental biquandle, see \cite{GEN} classify
virtual links up to mirror images? (We know that the biquandle has
the same value on the orientation reversed mirror image where the
mirror stands perpendicular to the plane (See \cite{HR,HRK}).

Are there good examples of weak biquandles which are not strong?

We would like to know more about the algebra with 2 generators $A,B$ and one
relation $[B,(A-1)(A,B)]=0$, (See \cite{FBu}). It is associated to the linear
case.

\item {\bf Virtualization and Unit Jones Polynomial.} Suppose the knot $K$ is classical and not trivial.
Suppose that ${\tilde K}$ (obtained from $K$ by virtualizing a subset of its crossings)
is not trivial and has a unit Jones polynomial , $V(\tilde K)=1$. Is it
possible that ${\tilde K}$ is classical (i.e. isotopic through virtual equivalence to a classical knot)?
\bigbreak

Suppose $K$ is a virtual knot diagram with unit Jones polynomial. Is $K$ equivalent to a classical diagram
via virtual equvalence plus crossing virtualization? 
(Recall that by crossing virtualization, we mean flanking a
classical crossing by two virtual crossings. This operation does not affect the value of the Jones polynomial.)
\bigbreak

Given two classical knots $K$ and $K',$
if $K$ can be obtained from $K'$ by a combination of crossing virtualization and
virtual Reidemeister moves, then is $K$ classically equivalent to $K'?$
\bigbreak

If the above two questions have affirmative answers, then the only classical knot with unit Jones polynomial is
the unknot.
\bigbreak

\item {\bf Virtual Quandle Homology.} Study virtual quandle homology in analogy to quandle homology 
\cite{CSH,GREENE}.

\item {\bf Khovanov Homology.} Construct a generalization of the Khovanov
complex for the case of virtual knots that will work for arbitrary virtual diagrams.
Investigate the Khovanov homology constructed in \cite{RBook, MAK}. The main construction in this approach 
uses an orientable atom condition to give a Khovanov homology over the integers for large classes of virtual links. 
The import of our question, is to investigate this structure and to possibly find a way to do Khovanov homology
for all virtual knots over the ring of integers. Similar questions can be raised for the presently evolving new classes of
Khovanov homolgy theories related to other quantum invariants.
\bigbreak

By a {\em K-full} 
virtual knot we mean a knot for which
there exists a diagram such that
the leading (the lowest, or both) term comes from the
$B$-state. Analogously, one defines the {\em Kho-full}
knot relative to the Khovanov invariant. Call such diagrams optimal diagrams.
(It is easy to find knots which are neither
K-full nor Kho-full.) 
\bigbreak

Classify all K-full (Kho-full)
knots.
\bigbreak

Are optimal diagrams always minimal with
respect to the number of classical crossings?
\bigbreak

Classify all diagram moves that preserve optimality. 
\bigbreak

Is it true that if a classical knot $K$
has minimal classical diagram with $n$ crossings
then any virtual diagram of $K$ has at least $n$
classical crossings?
\bigbreak

Can any virtual knot have torsion in
the $B$-state of the Khovanov homology (the genuine
leading term of some diagram)? Here we use the formulation of Khovanov
homology given in \cite{RBook} \cite{MAK}. 
\bigbreak

The behaviour of the lowest and the
leading term of the Kauffman bracket for virtual knots
was studied in \cite{RBook} and \cite{Avdeev} and \cite{NKamada}.
\bigbreak

\item {\bf Brauer algebra.}
The appropriate domain for the virtual recoupling theory is to place
the Jones-Wenzl projectors in the Brauer algebra. That is, when we add virtual crossings to the Temperley Lieb Algebra
to obrtain ``Virtual Temperley Lieb Algebra" the result is the Brauer algbra of all connections from $n$ points to 
$n$ points. What is the
structure of the projectors in this context? Can a useful algebraic generalization of the classical recoupling theory be formulated?

\item {\bf Virtual Alternating Knots.} Define and  classify
alternating virtual knots.

%\subitem 
Find an analogue of the Tait flyping conjecture and prove
it. Compare \cite{JZ}

%\subitem
Classify all alternating weaves on surfaces (without stabilization).

\item {\bf Crossing number problems.}

For each virtual link $L$, there are three crossing
numbers: the minimal number $C$ of classical crossings, the minimal
number $V$ of virtual crossings, and the minimal total number $T$ of
crossings for representatives of $L$. There are also a number of unknotting numbers: The
classical unknotting number is the number of crossing switches needed to unknot the knot (using any diagram for the knot).
The {\em virtual unknotting number} is the number of crossings one needs to convert from classical to virtual (by direct flattening)
in order to unknot the knot (using any virtual diagram for the knot). Very little is known.
Find out more about the virtual unknotting number.
\bigbreak

%\subitem 
What is the relationship between the least number of virtual
crossings and the least genus in a surface representation of the virtual knot.

%\subitem 
Is it true that $T=V+L$?

%\subitem 
Is there any algorithm for finding $V$ for some class of
virtual knots. For $T$, this is partially done for two classes of
links: quasialternating and some other, see \cite{Ma3}. For classical links and alternating diagrams see \cite{Mur,MT}.

%\subitem 
Are there some (non--trivial) upper and lower bounds for
$T,V,L$ coming from virtual knot polynomials?

\item {\bf Wild Virtuals.} Create the category of ``wild virtual knots''
and establish its axiomatics. In particular, one needs a theorem that states when
a wild equivalence of tame virtual links implies a tame equivalence of these links.

\item {\bf Vassiliev Invariants.} Understand the connection between virtual knot
polynomials and the Vassiliev knot invariants of virtual knots (in
Kauffman's sense). Some of that was done in \cite{VKT,GPV,SAW2,MAVI}.

The key question about this collection of invariants is this: {\em Does every Vassiliev invariant of finite type, for classical knots
extend to an invariant of finite type for long virtual knots?} Here we mean the problem in the sense of the formulation given in 
\cite{GPV}. In \cite{VKT} it was pointed out that there is a natural notion of Vassiliev invariants for virtual knots that has a different
notion of finite type from that given in \cite{GPV}. This alternate formulation needs further investigation.
\bigbreak

\item {\bf Embeddings of Surfaces.} Given a non-trivial virtual knot $K$. Prove that there
exists a minimal realization of $K$ in $N=S_{g}\times I$ and an
{\bf unknotted} embedding of $N\subset {\bf R}^{3}$ such that the
obtained classical knot in ${\bf R}^{3}$ is not trivial. (This problem is partially solved by 
Heather Dye in \cite{D1}.
\bigbreak

\item {\bf Non-Commutativity and Long Knots.} It is known that any classical long knot commutes with any
long knot. Is it true that it is the only case of commutativity
for virtual long knots. In other words, is it true that if $K$ and
$K'$ are long knots and $K{\#} K'$ is isotopic to $K'{\#} K$ then
there exists a virtual long knot $L$, classical long knots $Q,Q'$,
and non-negative integer numbers $m,n$ such that

$$K=L^{m} {\#} Q, \quad K'=L^{n} {\#} Q',$$

where by $L^{m}$ we mean the connected sum of $m$ copies of the
same knot.

\item {\bf The Rack Space.} The rack space was invented by Fenn,
Rourke and Sanderson \cite{FENN,FRS1,FRS2,FRS3}. The homology of the rack
space has been considered by the above authors and Carter, Kamada, Saito \cite{CSH}. For low dimensions,
the homology has the following interesting interpretations. Two
dimensional cycles are represented by virtual link diagrams
consistently colored by the rack, and three dimensional cycles by the
same but with the regions also colored.  See the Thesis of Michale
Greene \cite{GREENE}. So virtual links can give, in this way,
information about classical knots! For the second homology of the
dihedral rack, the results are given in Greene's Thesis. Computer
calculations suggest that for a prime $p$ the third homology has a
factor $Z_p$. Is this true in general?

Another line of enquiry is to look at properties of the birack space \cite{FRS2} and associated homology.

\end{enumerate}

\bigbreak

\bigbreak


\begin{thebibliography}{99}

\bibitem{P} 
M. M. Asaeda, J. H. Przytycki and A. S. Sikora, Categorification of the Kauffman bracket skein module of I-bundles over surfaces,
math.QA/0409414, {\em Algebraic and Geometric Topology} 4 (2004), paper no. 52, pp. 1177-1210.

\bibitem{Avdeev}
R.S.Avdeev, On extreme coefficients of the Jones-Kauffman
polynomial for virtual links, to appear in JKTR.

\bibitem{Bardakov}
V. G. Bardakov, The virtual and universal braids, arxiv:Math.GR/0407400 v1 23 Jul 2004.

\bibitem{FB}
A. Bartholomew and R. Fenn, Quaternionic invariants of virtual
knots and links, www.maths.sussex.ac.uk/Staff/RAF/Maths/Current/Andy/equivalent.ps,
 (preprint).

\bibitem{FBu} S.Budden and R.Fenn, The equation $[B,(A-1)(A,B)]=0$
and virtual knots and link, {\em Fund Math} 184 (2004) pp 19-29.


\bibitem{DB}
 D. Bar-Natan, On Khovanov's categorification of the Jones polynomial. {\em
Algebraic and Geometric Topology}, Vol. 2 (2002), pp. 337-370.

\bibitem{DBP}
 D. Bar-Natan, (private conversation).

\bibitem{Dror}
Dror Bar-Natan, Knot Atlas with Jeremey Green's atlas of virtual knots, http://www.math.toronto.edu/~drorbn/KAtlas/


\bibitem{CS}
 J.S.Carter and M. Saito, Diagrammatic invariants of knotted curves and surfaces,
(unpublished manuscript - 1992).

\bibitem{CSH}
J.S.Carter, S. Kamada, M. Saito, Geometric interpretations of quandle homology, {\em JKTR} {\bf 10}, No. 3 (2001), 345-386.

\bibitem{CSW1}
J.S.Carter, S. Kamada, M. Saito, Stable equivalence of knots on surfaces and virtual knot cobordisms,
math.GT/0008118, Knots 2000 Korea, Vol. 1 (Yongpyong). {\em J. Knot Theory Ramifications} 11 (2002), no. 3, 311--322.

\bibitem{Carter}
J.S. Carter, D. Jelsovsky, S. Kamada, M. Saito, Quandle homology groups,their Betti numbers and virtual knots,
arxiv:Math.GT/9909161 v1 28 Sep 1999.

\bibitem{P}
J. S. Carter and D. Silver (private conversation).

\bibitem{CT}
J. Conant and P. Teichner, Grope cobordism of
classical knots, (to appear).

\bibitem{DH}
 P. Dehornoy, ``Braids and Self-Distributivity", Progress in Math. vol. 192, Birkhauser
(2000).

\bibitem{Dye}
H. Dye, Characterizing Virtual Knots, Ph.D. Thesis (2002), UIC.

\bibitem{DKMin}
H. Dye and L. H. Kauffman, Minimal surface representations of virtual knots and links.
 {\em arXiv:math.AT/0401035 v1, Jan. 2004.} (to appear in Geometry and Topology)

\bibitem{DK}
 H. A. Dye and L. H. Kauffman, Virtual Knot Diagrams and the Witten-Reshetikhin-Turaev Invariant, math.GT/0407407,
(to appear in JKTR).

\bibitem{D1}
H. A. Dye, Non-Trivial Realizations of Virtual Link Diagrams, math.GT/0502477

\bibitem{D2}
H. A. Dye, Virtual knots undetected by 1 and 2-strand bracket polynomials, math.GT/0402308

\bibitem{EKT} S. Eliahou, L. Kauffman and M. Thistlethwaite, Infinite families of links with trivial Jones polynomial,
{\em Topology}, {\bf 42}, pp. 155--169.

\bibitem{FENN}
R. Fenn,
[www.maths.sussex.ac.uk//Staff/RAF/Maths/historyi.jpg], $(i = 1, 2, \cdots)$

\bibitem{FJK}
R. Fenn, M. Jordan and L. H. Kauffman, Biquandles and virtual links, {\em Topology and its Applications} {\bf 145} (2004), 157-175.

\bibitem{FR}
R. Fenn and C. Rourke, Racks and links in codimension two, {\em JKTR} No. 4, pp 343-406 (1992).

\bibitem{FRS1}
R. Fenn, C. Rourke and B. Sanderson, Trunks and classifying spaces, {\em Applied Categorical Structures}
3(1995) pp 321-356.

\bibitem{FRS2}
R. Fenn, C. Rourke and B. Sanderson, An introduction to species and the rack space, {\em Topics in Knot Theory},
Kluwer Acad. pp 33-55 (1993).

\bibitem{FRS3}
R. Fenn, C. Rourke and B. Sanderson, The rack space, arXiv:math.GT/0304228, (to appear in Trans. Amer. Math. Soc.).

\bibitem{FRR} 
R. Fenn, R. Rimanyi and C. Rourke, The Braid Permutation Group,
{\em Topology}, Vol. 36, No. 1 (1997), 123--135.

\bibitem{FT}
R. Fenn and P. Taylor, Introducing doodles, Lect. Notes in Maths. LMS No. 722, pp 37-43,

\bibitem{F}
R. Fenn, ``Techniques of Geometric Topology," LMS Lect. Notes
Series 57 (1983).

\bibitem{GQ}
R. Fenn (2005), Generalised Quaternions and Virtual Knots and Links, preprint.

\bibitem{GPV}
 Mikhail Goussarov, Michael Polyak and Oleg Viro
(2000), Finite type invariants of classical and virtual knots,
{\em Topology} {\bf 39}, pp. 1045--1068.


\bibitem{GREENE}
M. Greene, ``Some Results in Geometric Topology and Geometry", Thesis submitted for the degree of PhD, Warwick
Maths Institute, Sept. 1997.

\bibitem{HR}
D. Hrencecin, ``On Filamentations and Virtual Knot Invariants" Ph.D Thesis, Unviversity
of Illinois at Chicago (2001).

\bibitem{HRK} D. Hrencecin and L. H. Kauffman, ``On Filamentations and Virtual Knots", {\em Topology and Its Applications},
134(2003), 23-52.


\bibitem{JKS}
F. Jaeger, L. H. Kauffman and H. Saleur, The Conway polynomial in $R^{3}$
and in thickened surfaces: A new determinant formulation, {\em J. Comb. Theory Ser. B}
Vol. 61 (1994), 237-259.


\bibitem{KADOKAMI}
T. Kadokami, Detecting non-triviality of virtual links, {\em J. Knot Theory Ramifications} 12
(2003), no. 6, 781--803.

\bibitem{Kamada}
S. Kamada, Braid presentation of virtual knots and welded knots,
math.GT/0008092 (March 2000).

\bibitem{NKamada}
N. Kamada, Span of the Jones polynomial of an alternating virtual link,
{\em Alg and Geom Topology} Vol. 4 (2004) 1083-1101.

\bibitem{NKamamda1}
N. Kamada, On the Jones polynomials of checkerboard colorable virtual knots,
math.GT/0008074 

\bibitem{NSKamada}
N. Kamada and S. Kamada, Abstract link diagrams and virtual knots, {\em J. Knot Theory
Ramifications} 9 (2000), no. 1, 93--106.

\bibitem{KANENOBU} Kanenobu, T (2001), Forbidden moves unknot a virtual knot, {\em
Journal of Knot Theory and Its Ramifications}, {\bf 10} (1), pp.
89-96.


\bibitem{K}
L.H. Kauffman,
State Models and the Jones Polynomial,
{\em Topology} {\bf 26} (1987), 395--407.

\bibitem{KS}
L. H. Kauffman and H. Saleur, Free fermions and the Alexander-Conway polynomial, {\em Comm. Math.
Phys.} 141, 293-327 (1991).

\bibitem{GAUSS}
L. H. Kauffman, Gauss Codes, quantum groups
and ribbon Hopf algebras, {\em Reviews in Mathematical Physics}
{\bf 5} (1993), 735-773. (Reprinted in \cite{KNOTS}, 551--596.


\bibitem{KNOTS}
L. H. Kauffman, ``Knots and Physics", World Scientific,
Singapore/New Jersey/London/Hong Kong, 1991, 1994, 2001.

\bibitem{TL}
 L. H. Kauffman and S. Lins, ``Temperley - Lieb Recoupling Theory and
Invariants of 3-Manifolds".  Princeton University Press - 1994.

\bibitem{KLOM}
 L. H. Kauffman and S. J. Lomonaco, Quantum entanglement and topological entanglement. New J. Phys. 4 (2002), 73.1 - 73.18.

\bibitem{IM}
 L. H. Kauffman and D. E. Radford, Invariants of 3-manifolds derived from finite
dimensional Hopf algebras.  Journal of Knot Theory and its Ramifications, Vol.4, No. 1 (1995), pp.
131-162.

\bibitem{KFI}
 L. H. Kauffman, Vassiliev invariants and functional integration without integration. {\em Stochastic
 analysis and mathematical physics} (SAMP/ANESTOC 2002), 91--114, World Sci. Publishing, River Edge, NJ, 2004..

\bibitem{KL}
L. H. Kauffman and S. Lambropoulou, Virtual Braids, math.GT/0407349 (to appear in Fund. Mathematica)

\bibitem{KaMa}
 L. H. Kauffman, V.O. Manturov, Virtual biquandles,
 math.GT/0411243(to appear in Fund. Mathematica)


\bibitem{RI}
 L. H. Kauffman,  Right Integrals and Invariants of Three-Manifolds, Proceedings of Conference in
Honor of Robion Kirby's 60th Birthday, Geometry and Topology Monographs, Vol. 2 (1999),
215-232.

\bibitem{KRO}
 L. H. Kauffman and David E. Radford, Oriented quantum algebras and
invariants of knots and links,  {\em Journal of Algebra}, Vol. 246, 253-291 (2001).

\bibitem{KRCAT}
 L. H. Kauffman and David E. Radford, Oriented quantum algebras,
categories and invariants of knots and links.  {\em JKTR}, vol 10, No. 7 (2001), 1047-1084. \smallbreak

\bibitem{GEN}
 L. H. Kauffman and D. E. Radford, Bi-oriented Quantum Algebras, and a
Generalized Alexander Polynomial for Virtual Links, {\em Diagrammatic Morphisms and Applications} (San Francisco, CA, 2000), 113--140,
 Contemp. Math., 318, Amer. Math. Soc., Providence, RI, 2003..


\bibitem{KDL}
L. H. Kauffman,D. De Wit and J. Links. On the
Links-Gould Invariant of Links. JKTR, Vol. 8, No. 2 (1999), 165-199.

\bibitem{INT}
L. H. Kauffman, Knot theory and the heuristics of functional integration.
 Physica A 281 (2000) 173-200.

\bibitem{VKT}
L. H. Kauffman, Virtual Knot Theory , {\em European J. Comb.} (1999) Vol. 20, 663-690.


\bibitem{SVKT}
L. H. Kauffman, A Survey of Virtual Knot Theory in {\em Proceedings of Knots in Hellas '98},
World Sci. Pub. 2000 , pp. 143-202.

\bibitem{DVK}
L. H. Kauffman, Detecting Virtual Knots, in Atti. Sem. Mat. Fis.
Univ. Modena Supplemento al Vol. IL, 241-282 (2001).

\bibitem{SLK}
L. H. Kauffman, A Self-Linking Invariant of Virtual Knots
math.GT/0405049 (to appear in Fund. Mathematica)

\bibitem{KD}
L. H. Kauffman, Knot Diagrammatics, math.GN/0410329 (to appear in Handbook of Knot Theory)

\bibitem{DKT}
L. H. Kauffman, Diagrammatic Knot Theory (in preparation).

\bibitem{LF} L. H. Kauffman, Long Flat Virtual Knots (in preparation).

\bibitem{KM}  L. H. Kauffman and V. O. Manturov, Virtual Knots and Links, (To Appear
in Proceedings of the Steklov Inst. RAS).

\bibitem{KHO}
M. Khovanov, Doodle groups, {\em Trans. Amer. Math. Soc.}, Vol. 349, No. 6, pp 2297-2315.


\bibitem{KiSa} 
T. Kishino and S. Satoh, A note on non-classical virtual knots,
{\em J. Knot Theory Ramifications} 13 (2004), no. 7, 845--856.

\bibitem{KIS}
T. Kishino, 6 kouten ika no kasou musubime no
bunrui ni tsiuti (On classification of virtual links whose crossing
number is less than or equal to 6), Master Thesis, Osaka City
University, 2000.

\bibitem{knotilus} Knotilus website: http://srankin.math.uwo.ca/cgi-bin/retrieve.cgi/html/start.html

\bibitem{KUP}
Greg Kuperberg, What is a virtual link? arXiv:math.GT $\slash$
0208039 v1 5 Aug 2002,{\em Algebraic and Geometric Topology},
2003, {\bf 3}, 587-591.

\bibitem{LM}
 W.B.R. Lickorish and K.C. Millett,
Some evaluations of link polynomials, {\em Comment. Math. Helvetici} {\bf 61} (1986),
349--359.


\bibitem{MALONG}
 V. O. Manturov, Long virtual knots and their
invariants, {\em J. Knot Theory Ramifications} 13 (2004), no. 8,
1029--1039.

\bibitem{Ma0}
V. O. Manturov, (2004), Knot Theory, CRC-Press,
Chapman\& Hall.

\bibitem{Ma1}
V. O. Manturov, (2002), Invariants of Virtual
Links, {\em Doklady Mathematics} Vol. 65 (3),2002, 329-332.

\bibitem{MACURVES'}
V.O. Manturov, (2003), Curves on
Surfaces, Virtual Knots, and the Jones--Kauffman Polynomial, {\em
Doklady Mathematics} Vol. 65 (3), 2003, 326-329.  

\bibitem{MAKH}
V.O. Manturov (2004), The Khovanov polynomial for virtual
knots, {\em Doklady Mathematics} Vol. 70, No. 2, 2004, 679-682. 

\bibitem{MAK} 
V. O. Manturov, Khovanov Complex for
Virtual Links, Arxiv/Math:GT/0501317 (to appear in 
Fundamental and Applied Mathematics, 2005).

\bibitem{Ma3}
V. O. Manturov, (2003), Atoms and minimal diagrams
of virtual links, {\em Doklady Mathematics}, {\bf 391} (2),
pp.136-138

\bibitem{MAPOLY}
V.O. Manturov,(2002), Two--variable invariant
polynomial for virtual links, {\em Russian Math. Surveys},{\bf 57},
(5), pp. 997-998.

\bibitem{MAPOLY'}
V. O. Manturov, (2003), Multivariable polynomial
invariants for virtual knots and links, {\em Journal of Knot
Theory and Its Ramifications}, {\bf 12},(8), pp. 1131-1144


\bibitem{MACURVES}
V. O. Manturov, (2003), Kauffman--like
polynomial and curves in $2$--surfaces, {\em Journal of Knot
Theory and Its Ramifications}, {\bf 12},(8),pp.1145-1153.


\bibitem{MAVI}
V.O. Manturov, Vassiliev Invariants for Virtual
Links, Curves on Surfaces, and the Jones-Kauffman polynomial, to
appear in JKTR


\bibitem{Ma12} 
Finite-type invariants of Virtual Links and the
Jones-Kauffman polynomial, {\em Doklady Mathematics} Vol. 69 (2), 2004, pp. 164-167.

\bibitem{Ma7}
V. O. Manturov, (2002), On Invariants of Virtual
Links, {\em Acta Applicandae Mathematicae}, {\bf 72} (3), pp.
295-- 309.

\bibitem{MaLong} V.O.Manturov, On long virtual knots, Doklady
Mathematics, 2005



\bibitem{Ma9}
V.O. Manturov (2003)., O raspoznavanii virtual'nyh
kos (On the recognition of virtual braids), POMI Scientific
Seminars. Geometry and Topology.8. Saint-Petersburg, pp.
267-286.


\bibitem{Ma10}
 V.O. Manturov, (2003), Atoms and Minimal Diagrams
of Virtual Links, Doklady Mathematics, {\bf 68}, (1), pp.
37-39.

\bibitem{Ma11}
V.O.Manturov (2004) Flat Hierarchy, (to appear in Fundamenta Mathematicae)

\bibitem{MAN}
V. O. Manturov, Compact and long virtual knots, to appear in
{\em Proceedings of the Moscow Mathematical Society}.

\bibitem{RBook} V.O.Manturov, Teoriya Uzlov (Knot Theory),
RCD, Moscow-Izhevsk, 2005, 512pp, (in Russian).


\bibitem{MT} W. Menasco and M.B. Thistlethwaite,
The classification of alternating links,
{\em Ann. of Math.} {\bf 138} (1993), 113--171.


\bibitem{Mur}
K. Murasugi,
The Jones polynomial and classical conjectures in knot theory,
{\em Topology} {\bf 26} (1987), 187--194.

\bibitem{NELSON}
S. Nelson, Unknotting virtual knots with Gauss diagram forbidden moves. {\em JKTR} 10 (2001), no. 6,931-935.

\bibitem{Nelson1}
S. Nelson,Virtual crossing realization, math.GT/0303077

\bibitem{SATOH}
S. Satoh, Virtual knot presentation of ribbon
torus-knots, JKTR, Vol. 9 No. 4 (2000), pp. 531-542.

\bibitem{SAW}
 J. Sawollek, On Alexander-Conway polynomials for virtual knots and links,
arXiv:math.GT/9912173, 21 Dec 1999.


\bibitem{SAW2}
 J. Sawollek, An orientation-sensitive Vassiliev
invariant for virtual knots, 2002, arXiv:math.GT/ 0203123 v3.

\bibitem{Schellhorn}
W. J. Schellhorn, Filamentations for virtual links, 
math.GT/0402162
 

\bibitem{SW}
 D. S. Silver and S. G. Williams, Alexander Groups and Virtual Links, {JKTR}, vol. 10,
(2001), 151-160.

\bibitem{SW1}
 D. S. Silver and S. G. Williams, Alexander groups of long virtual knots
math.GT/0405460

\bibitem{SW2}
D. S. Silver and S. G. Williams, On a class of virtual knots with unit Jones polynomial,
{\em J. Knot Theory Ramifications} 13 (2004), no. 3, 367--371.

\bibitem{SW3}
D. S. Silver and S. G. Williams, Polynomial invariants of virtual links. J. Knot Theory
Ramifications 12 (2003), no. 7, 987--1000.

\bibitem{MT}
 M. Thistlethwaite, Links with trivial Jones polynomial. {JKTR} 10 (2001), no. 4, 641--643.

\bibitem{T}
M.B. Thistlethwaite,
A spanning tree expansion of the Jones polynomial,
{\em Topology} {\bf 26} no. 3 (1987), 297--309.

\bibitem{TURAEV}
 V. Turaev, Virtual strings and their cobordisms,
math.GT/0311185.

\bibitem{Ver}
V. V. Vershinin, On homology of virtual braids and Burau representation, 
math.GT/9904089

\bibitem{W}
S. Winker. PhD. Thesis, University of Illinois at
Chicago (1984).

\bibitem{WITT}
 E. Witten. Quantum Field Theory and the Jones Polynomial. Comm. in Math. Phys.
Vol. 121 (1989), 351-399. 

\bibitem{JZ}
P. Zinn-Justin and  J. B. Zuber, Matrix integrals and the generation and counting of virtual tangles and links,
{\em J. Knot Theory Ramifications} 13 (2004), no. 3, 325--355.

\bibitem{JZ1}
P. Zinn-Justin and  J. B. Zuber, Tables of Alternating Virtual Knots - http://ipnweb.in2p3.fr/~lptms/membres/pzinn/virtlinks/

\end{thebibliography}
\end{document}